\documentclass[11pt]{amsart} 
\usepackage{amsfonts}
\usepackage{amsmath}
\usepackage{amssymb}
\usepackage{amsthm}
\usepackage{cite}
\usepackage{epsfig} 
\usepackage{rotate}
\usepackage{tikz}
\usepackage{pgfplots}
\usetikzlibrary{math}

\newtheorem{theorem}{Theorem} 
\newtheorem{corollary}[theorem]{Corollary}
\newtheorem{lemma}[theorem]{Lemma}
\newtheorem{proposition}[theorem]{Proposition}

\def\R{\mathbb R}
\def\Q{\mathbb Q}
\def\Z{\mathbb{Z}}

\def\a{\mathrm{a}} 
\def\b{\mathrm{b}}
\def\c{\mathrm{c}} 
\def\cA{\mathcal{A}}
\def\cC{\mathcal{C}}
\def\Cbar{\tilde C}
\def\bC{\mathbf{C}}
\def\bc{\mathbf{c}}
\def\be{\mathbf{e}}
\def\cD{\mathcal{D}}
\def\rf{\mathrm{f}}
\def\rF{\mathrm{F}}
\def\cG{\mathcal{G}}
\def\cL{\mathcal{L}}
\def\rM{\mathrm{M}}
\def\cQ{\mathcal{Q}}
\def\rR{\mathrm{R}}
\def\bv{\mathbf{c}}

\def\bd{\mathrm{dom}}

\def\cp{{U}}
\def\idl#1{\langle #1\rangle}
\def\n2{\mathbf{k}}

\def\qad{\hskip 5pt}
\def\mod#1{\,({\rm mod\ }#1) }
\def\proof{\noindent {\sc Proof.}\hskip 10pt}
\def\endproof{\hfill\vbox{\hrule \hbox{\vrule\kern4pt\vbox{\kern4pt
\kern4pt}\kern4pt\vrule}\hrule}\smallskip}

\begin{document}
\title[]{Critical curves of a piecewise linear map}
\author{John A G Roberts} 
\address{School of Mathematics and Statistics,
University of New South Wales,
Sydney, NSW 2052, Australia}
\email{jag.roberts@unsw.edu.au}
%
\author{Asaki Saito}
\address{Future University Hakodate,
116--2 Kamedanakano-cho,
Hakodate, Hokkaido 041--8655, Japan}
\email{saito@fun.ac.jp}
\author{Franco Vivaldi}
\address{School of Mathematical Sciences, Queen Mary,
University of London,
London E1 4NS, UK}
\email{f.vivaldi@maths.qmul.ac.uk}
%
\begin{abstract}
We study the parameter space of a family of planar maps, 
which are linear on each of the right and left half-planes.
We consider the set of parameters for which every orbit recurs 
to the boundary between half-planes. These parameters consist of
algebraic curves, determined by the symbolic dynamics of the itinerary 
that connects boundary points. 
We study the algebraic and geometrical properties of these curves, in
relation with such a symbolic dynamics.
\end{abstract}
\date{\today}
\maketitle

\section{Introduction}\label{section:Introduction}

Piecewise linear (affine) maps are maps defined on a partitioned phase space where a different linear (affine) map acts on each region of the partition. Often this can be achieved to maintain continuity of the map. Dissipative versions arise naturally in engineering and physical models and the nature of the attractors in them and the possible bifurcations have received considerable attention \cite{diBernardoEtAl}.  

Less well-studied is the conservative case; in particular, regular 
orbits in piecewise linear symplectic maps are not well understood.
A much studied two-parameter family on the plane is 
\cite{Devaney,BeardonBullettRippon,LagariasRainsI,
LagariasRainsII,LagariasRainsIII,Garcia-MoratoEtAl}
\begin{equation}\label{eq:LR}
\rF(x,y)=\begin{cases} 
(\a x-y,x)& \,\,x>0\,\lor\,(x=0 \land y\leqslant 0)\\
(\b x-y,x)& \,\,\mbox{otherwise}
\end{cases}
\end{equation}
where $\a$ and $\b$ are real parameters.
(Our definition of $\rF$ on the line $x=0$ is a variant of that found
in the literature.) The map $\rF$ can lay claim to being the normal form for a piecewise linear map acting on the partition into the left and right half-planes \cite{Garcia-MoratoEtAl}.

\begin{figure}[t]
  \begin{minipage}{5cm}
          \centering
         \includegraphics[scale=0.21]{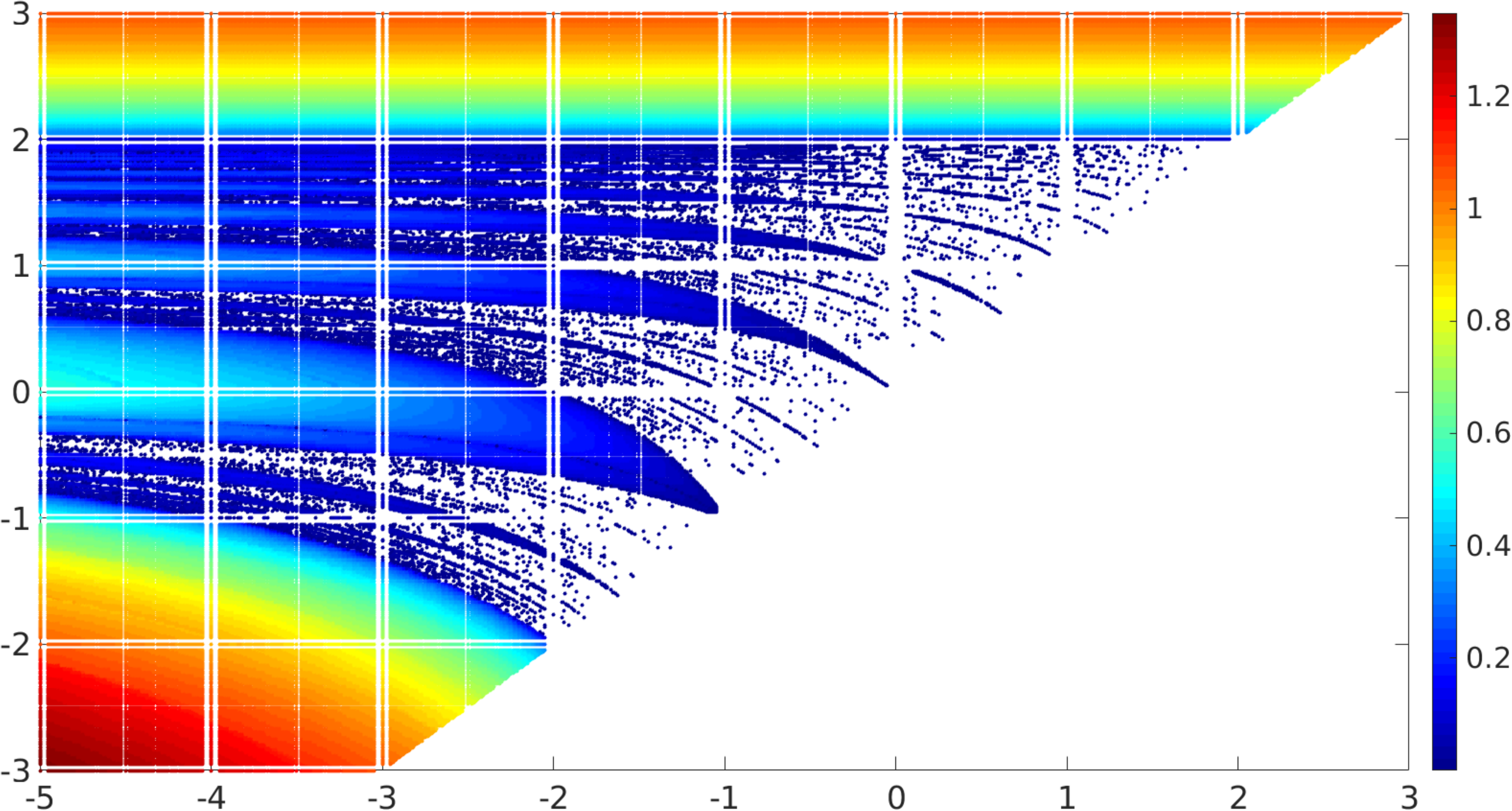}\\
  \end{minipage}
\hspace*{100pt}
  \begin{minipage}{5cm}
          \centering
          \includegraphics[scale=0.30]{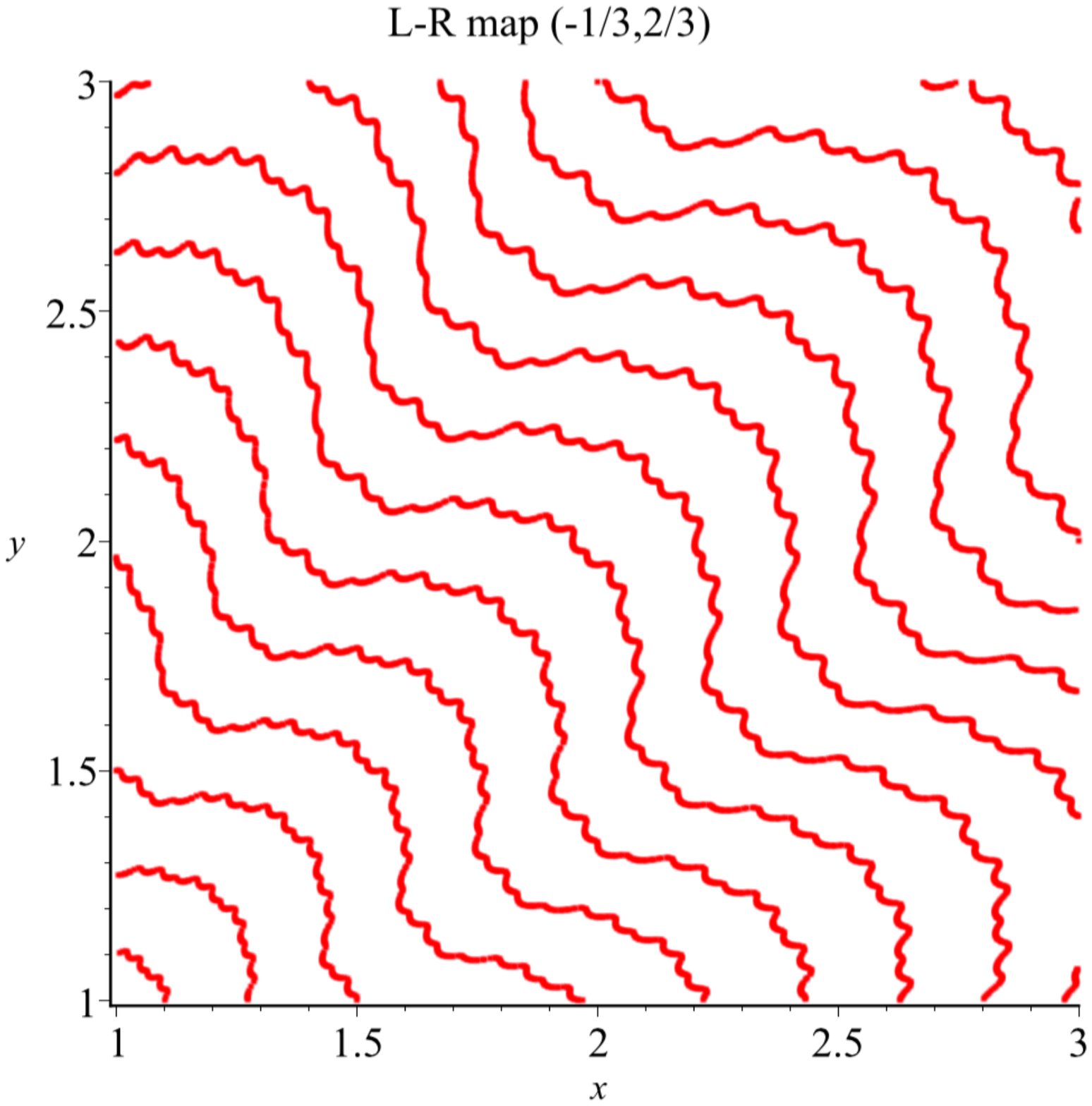}\\
  \end{minipage}
\hspace*{5pt}
\vspace*{-50pt}
\hfil
\caption{\label{fig:LRoverview}\rm\small
Left: the parameter space of the map (\ref{eq:LR}), with 
mode-locking regions colour-coded by the escape rate (only 
the range $\a\leqslant \b$ is shown).
Right: some orbits for the parameter pair $(\a,\b)=(-1/3,2/3)$, 
seemingly conjugate to an irrational rotation with 
symbolic complexity $\mathcal{K}(n)=2n$.
}
\end{figure}

The map $\rF$ sends rays through the origin into themselves, and the
ray dynamics is a smooth circle map $\rf$. The purpose of this 
work is to study the symbolic dynamics of $\rf$ with respect to 
the above binary partition of the plane,
extending the works \cite{LagariasRainsI,LagariasRainsII} on the system (\ref{eq:LR}), 
and complementing the works \cite{SimpsonMeiss,Simpson:17,Simpson:18,Garcia-MoratoEtAl}
on mode-locking and bifurcations in piecewise-linear maps.
This work is the first part of a planned broader study \cite{RobertsEtAl}.

Each orbit of $\rF$ corresponds to developing a product or {\em word} in the matrices  
$A=\left(\begin{smallmatrix}\a &-1\\ 1&0\end{smallmatrix}\right)$ and
$B=\left(\begin{smallmatrix}\b &-1\\ 1&0\end{smallmatrix}\right)$, 
following the symbol sequence of the dynamics between the left and 
right-half planes. 
(We use roman fonts for the real parameters $\a,\b$ and italic fonts for 
the symbols $a,b$ of the corresponding half-planes; the latter appear
as letters in words or indeterminates in algebraic expressions.)
Relevant observables are the frequencies of the two type of matrices in the developing word, which are orbit-dependent.  
In this sense, $\rF$ generalises a model in condensed matter physics,  the discrete Schr\"odinger equation on 
a one-dimensional lattice  (see, e.g., \cite{Roberts} 
and references therein).  
Here one propagates the solution along the lattice by developing a matrix word in 
the aforementioned matrices $A$ and $B$ (with $a$ and $b$  being physically relevant 
parameters) according to a two-letter substitution rule which prescribes the 
asymptotic frequencies of $A$ and $B$ in the word.  
The case of quasiperiodic words (e.g., the Fibonacci sequence) has received much attention. 

Matrix words of the type described also arise naturally in the classical study of continued fraction expansions in number theory where $\a$ and $\b$ are taken to be integers \cite{MGO}. Nevertheless, the theory of {\em continuant polynomials} which is used to study continued fractions can equally well be applied to study the dynamics of $\rF$ with $\a, \b \in \R$ \cite{B-MEtAl, GrahamEtAl, MGO}.

Our knowledge of parameter space of $\rF$ is limited to four countable families
of algebraic curves $C(a,b)=0$, where $C$ is a polynomial with integer 
coefficients. The first two families are known explicitly, 
while the other two are constructed using certain finite symbolic codes, 
whose general form is still unavailable. They are:
\begin{enumerate}
\item [i)]  An infinite family of line segments (together with their 
image under parameter exchange\footnote{Exchanging $\a$ and $\b$ in 
(\ref{eq:LR}) leads to a conjugate system, so it suffice to study 
the case $\a\leqslant \b$.}), where the rotation number of $\rf$ 
is known \cite{BeardonBullettRippon,LagariasRainsI}.
\item [ii)] An infinite family of hyperbolae, invariant under parameter 
exchange, where the map $\rf$ is of finite order and the rotation number 
is constant 
\cite{BeardonBullettRippon,LagariasRainsI}.
\item [iii)] A family of algebraic curves defined by a dynamical condition,
for which the map $\rF$ is known to support invariant curves with rational or 
irrational rotation number, consisting of finitely many arcs of conic 
sections glued together \cite{LagariasRainsII}.
\item [iv)] A family of algebraic curves defined by a dynamical 
condition, bounding the mode-locking regions
\cite{SimpsonMeiss,Simpson:17,Simpson:18,Garcia-MoratoEtAl}.
\end{enumerate}

A central question regards the existence and properties of 
quasi-periodic orbits whose closure are topological circles, 
special cases of which occur for parameters of type i) and iii).
General results are scarce.
M. Herman's work in the broader setting of the Froeschl\'e group
\cite[Theorem VIII.5.1]{Herman} implies that any map $\rF$ having 
an irrational rotation number $\theta=\theta(\a,\b)$ with bounded
partial quotients in its continued fraction expansion is topologically 
conjugate to a rotation of the plane, and hence has invariant circles.
This set of parameter values has zero two-dimensional Lebesgue measure. 
The results of \cite{LagariasRainsII} on curves of type iii) are 
complementary, as they regard one-dimensional sets of parameters (see below).
Numerical experiments suggest ubiquity of parameter pairs for which all 
orbits are dense on non-smooth curves, as in figure \ref{fig:LRoverview}, 
right.

The present work is devoted to the development of a theory of the
curves of type iii), which we call the \textbf{critical curves} 
in parameter space.
They are defined by the
condition that an initial boundary ray (the positive or negative ordinate 
semi-axis in $\R^2$) be sent by the map $\rF$ to the opposite boundary ray in 
a prescribed number of iterations. If the parameters $\a,\b$ belong to critical curve,
and if the rotation number of the corresponding circle map $\rf$ is irrational, 
then the phase space of $\rF$ foliates into piecewise smooth invariant curves, 
consisting of arcs of conic sections joined together 
\cite[theorem 2.2]{LagariasRainsII}\footnote{The authors of \cite{LagariasRainsII}
do not introduce explicitly the notion of critical curve.}. 

We now summarise the contents and main results of this paper.
In section \ref{section:Background} we provide some background on 
the symbolic dynamics of rotations and on the circle map $\rf$. We then give some
preliminary results for the dynamics of $\rF$ on a half-plane and its relation to the parameter.
In the following section 
we introduce the main object of
study, the critical curves, which correspond to the occurrence of Sturmian-type 
words in the symbolic dynamics (see \cite[section 6]{PytheasFogg}).
Using the theory of continuants, we 
determine some general algebraic properties of these curves
(propositions \ref{prop:ContinuedFractions},\ref{prop:Palindrome}, 
section \ref{section:AlgebraicProperties}).
We link the symbolic dynamics to the time-reversal 
symmetry of the map $\rF$ (theorem \ref{thm:PalindromeII}),
and show that critical curves have several disjoint branches,
whose number is determined by the code. We establish that
only one branch ---indeed, only a part of it--- is relevant 
to the dynamics of the map $\rF$ (theorems \ref{thm:Branches} and \ref{thm:LegalEagle}),
meaning that the code that defines the curve is the symbolic 
dynamics of a segment of an actual orbit of the map $\rF$, 
with initial condition of the appropriate boundary ray.
Subsequently (section \ref{section:Functions})
we consider functions defined over critical curves, and 
derive formulae for the Poisson brackets of two curves 
(theorem \ref{thm:Poisson}), for later use.

In section \ref{section:Intersections} we consider intersections of 
critical curves, which lead to periodic orbits. We introduce the 
concept of intersection sequence, with which we classify the 
intersections of a curve with curves of lower degree generated 
by sub-words (theorem \ref{thm:DoublePoints}).
We construct an associated geometrical object, called the polygonal 
of the curve, by means of which we formulate sufficient conditions 
for the transversality 
of these intersections (theorem \ref{thm:Regular}). 
We also formulate conditions under which intersections
delimit the portion of a curve which has dynamical significance.

In section \ref{section:RotationalDomains} we consider the
first generation of curves of type iii). These are obtained by
concatenating the symbolic words of two curves of type i),
and then allowing for repetitions of the concatenated word (theorem
\ref{thm:FirstGenerationCurves}).

Identities of Chebyshev polynomials are collected in an appendix.

\maketitle

\section{Background}\label{section:Background}

\subsection{Rotational words}\label{section:RotationalWords}
For references for this section, see, e.g., \cite{PytheasFogg,AlessandriBerthe}.

The \textbf{rotational words} are the symbolic dynamics of 
rotations with respect to a two-element partition of the circle.
A rotational orbit has the form $x_t=\{\theta t+x_0\}$
for some $\theta,x_0\in [0,1)$, where $\{\cdot\}$ denotes
the fractional part.
Without loss of generality, we choose the partition $I_a=[0,\rho)$
and $I_b=[\rho,1)$, the subscript denoting the symbol associated
to each interval.
Thus a rotational word is determined by a triple $(\theta,\rho,x_0)$.
Rotations are invertible, so all rotational words can be extended 
to the left; the extension is unique, apart from some notable 
special cases (see below). 

If we denote by $\cL_n$ the collection of all sub-words 
(or \textit{factors}) of length $n$ in an infinite word $W$,
then the \textit{complexity function} of $W$ is $\mathcal{K}(n)=|\cL_n|$.
An infinite word is eventually periodic if and only if 
$\mathcal{K}(n)\leqslant n$ for some $n$, so we require 
$\mathcal{K}(n)\geqslant n+1$. 
For rotational words, $\mathcal{K}(n)$ is (eventually) an affine 
function of $n$. The precise form of $\mathcal{K}$ depends on 
whether or not the boundary points $0$ and $\rho$ of the 
partition are on the same doubly-infinite orbit.
If they are, we let $\ell$ be the transit time
from one to the other; if they aren't, we let $\ell=\infty$.
Then, irrespective of the initial condition $x_0$, the rotational
word has complexity \cite[theorem 10]{AlessandriBerthe}
\begin{equation}\label{eq:K(n)}
\mathcal{K}(n)=\begin{cases} 2n& n\leqslant \ell \\n+\ell & n>\ell.\end{cases}
\end{equation}
The case $\ell=1$ are the much studied \textit{Sturmian} words \cite[section 6]{PytheasFogg}.
If $1<\ell<\infty$, we call the word \textit{quasi-Sturmian} (of length $\ell$), and
\textit{free} if $\ell=\infty$.

All irrational rotations are minimal (and indeed uniquely ergodic), and
this implies that for fixed irrational $\theta$ and arbitrary $\rho$,
each factor of any rotational word occurs infinitely often and with
bounded time between successive occurrences, hence with a well-defined 
frequency; the latter may be computed as the length of an 
interval \cite{AlessandriBerthe}.
The symbolic language [the sequence $(\cL_n)$] does not 
depend on the initial condition $x_0$ \cite[p.~105]{PytheasFogg}, 
and for our purpose the parameter space 
for rotations will be the unit square $[0,1)^2$.
The parameters $\theta$ and $\rho$ will be referred to as the 
\textbf{rotational parameters}. 

If the points $0$ and $\rho$ are on the same orbit, separated by $\ell$ 
iterations ($\ell$ could be negative),
then $\ell\theta\equiv\rho\mod{1}$, that is, the quasi-Sturmian 
configurations for a given $\ell$ correspond to a finite collection of 
segments in rotational parameter space.

\subsection{The circle map}\label{section:CircleMap}

Let $\rf:[0,1)\to[0,1)$ be the circle map associated to (\ref{eq:LR}), namely
\begin{equation}\label{eq:CircleMap}
\rf(e^{2\pi i x})=\frac{\rF(e^{2\pi i x})}{\Vert \rF(e^{2\pi i x})\Vert}
\qquad e^{2\pi i x}\cong (\cos(2\pi x),\sin(2\pi x))\in\mathbb{R}^2.
\end{equation}
This is an orientation-preserving homeomorphism whose derivative is
continuous and of bounded variation \cite[theorem 3.1]{LagariasRainsI}.
Hence $\rf$ has a well-defined rotation number $\theta=\theta(\a,\b)$. The latter is
a continuous function of the parameters, and its range is the interval 
$[0,1/2]$ \cite[theorems 2.1 and  3.3]{LagariasRainsI}.
From Denjoy theorem \cite[p.~401]{KatokHasselblatt}, if $\theta$ is irrational
then $\rf$ is topologically conjugate to a rotation by $\theta$, which in turn
determines the symbolic language $\cL$, irrespective of the initial conditions.

We divide the domain of the circle map $\rf$ into two intervals $I_a$ and $I_b$, 
and associate to an orbit a symbolic word $W$ in the letters $a$ and $b$, 
determined by the two branches of the map $\rF$ in (\ref{eq:LR}). 
We denote by $W_n$ the prefix of $W$ of length $n$, and by $|W_n|_w$ 
the number of times the finite word (factor) $w$ appears in $W_n$.
The rotation number $\theta(\a,\b)$ and the density $\rho(\a,\b)$ are
defined as follows:
\begin{equation}\label{eq:ThetaRho}
\theta(\a,\b)=\lim_{n\to\infty}\frac{1}{n}|W_n(\a,\b)|_{ab},
\qquad
\rho(\a,\b)=\lim_{n\to\infty}\frac{1}{2n}\bigl(|W_n^+(\a,\b)|_a+|W_n^-(\a,\b)|_a\bigr),
\end{equation}
where $W$ is arbitrary and $W^\pm$ is the word of the orbit with 
initial condition $(0,\pm1)$.

The above limits exist because they are frequencies
of factors of rotational words.
Since every appearance of the factor $ab$ corresponds to one loop around 
the origin, our definition of rotation number is equivalent to the 
standard one for circle homeomorphisms \cite[chapter 11]{KatokHasselblatt};
in particular, the rotation number is independent from the initial condition. 
This also holds for the density, as long as the rotation number is irrational,
from minimality of irrational rotations (section \ref{section:RotationalWords}).
If the rotation number is rational, then the density may depend on the initial
condition; in our definition (\ref{eq:ThetaRho}) the union of the orbits corresponding
to $W^\pm$ is invariant under time-reversal symmetry [because $\rF(0,\mp1)=(\pm1,0)$].
This ensures that $\rho(\a,\a)=1/2$ for all real $\a$, so that $\rho$ is
continuous on the main diagonal in parameter space.

We consider the level sets of the rotation number:
\begin{equation}\label{eq:LevelSets}
\begin{array}{rcl}
\Theta(x)&=&\{(\a,\b)\in\mathbb{R}^2\,|\, \theta(\a,\b)=x\}
\end{array}
\qquad 0\leqslant x \leqslant 1/2.
\end{equation}
If $x=p/q$ is rational, then $\Theta(p/q)$ will be called a 
\textit{resonance}.\footnote{We prefer this term to the customary
terms \textit{tongue} or \textit{sausage}.}
A point $(\a,\b)$ of a resonance is a \textit{pinch-point} if the set
$\Theta(p/q)\setminus \{(\a,\b)\}$ is locally disconnected near $(\a,\b)$.
Note that the density $\rho$ is not necessarily continuous within a resonance,
because the same rotation number may be associated to more than one code.
For instance, the periodic codes $aabb$ and $abbb$ have the same rotation number 
1/4, but distinct densities 1/2 and 1/4.

We now remove from the parameter space $\mathbb{R}^2$ of the circle map 
$\rf$ the interior of the resonances $\Theta(0)$ and $\Theta(1/2)$, given by
\begin{equation}\label{eq:Resonances}
\Theta(1/2)=\{(\a,\b)\in\mathbb{R}^2\,|\,\a\b\leqslant 4, \a,\b<0\},\qquad
\Theta(0)=\{(\a,\b)\in\mathbb{R}^2\,|\,\max(\a,\b)\geqslant 2\}.
\end{equation}
We obtain an infinite strip $\cA$,
the region lying between the boundaries of $\Theta(0)$ and $\Theta(1/2)$.
(If we identify the points $(-\infty,\zeta)$ and $(\zeta,-\infty)$, 
for $0\leqslant\zeta\leqslant 2$, then $\cA$ becomes a compact set
---a topological annulus.)

A detailed investigation of parameter space outside resonances
requires replacing original parameters $\a,\b$ with the
rotational parameters $\theta,\rho$, which are suited for
an arithmetical analysis of critical curves. 
This is the subject of a forthcoming investigation \cite{RobertsEtAl}.
\subsection{Basic dynamical properties}\label{section:Legal}

The area-preserving map $\rF$ of \eqref{eq:LR} has a common form of 
a \textbf{reversible} map, i.e., a map conjugate 
to its inverse via an involution, specifically
\begin{equation}\label{eq:defR}
\rF^{-1} = \rR \circ \rF \circ \rR, \quad \rR: (x,y)\mapsto (y,x).
\end{equation}
Reversible maps are well-studied \cite{LambRoberts} and an equivalent definition is that they can be written as the composition of two involutions (e.g., $\rF$ is the composition of $\rF \circ \rR$ and $\rR$). An orbit of $\rF$ is called {\em symmetric} if it is $\rR$-invariant and {\em asymmetric} otherwise, in which case it forms an asymmetric pair with its $\rR$-image. A symmetric orbit must contain one or two points from the {\em symmetry lines}
$$ Fix(\rR)  \cup  Fix(\rF \circ \rR),$$
and is periodic if and only if it contains two points, where
\begin{equation}\label{eq:defFix}
Fix(\rR):=\{(x,x): x \in \R \},\quad Fix(\rF \circ \rR):=\{(\frac{\a}{2}\,y,y): y > 0\} \cup \{(\frac{\b}{2}\,y,y): y \leqslant 0 \}.
\end{equation}
 
As also realised in \cite{LagariasRainsI}, a special role in the dynamics of $\rF$ is played by the positive and negative ordinate semi-axes and we define:
\begin{equation}\label{eq:defLpm}
L^{+} := \{(0,y): y \geqslant 0\}, \quad
L^{-} := \{(0, -y): y \geqslant 0\}.
\end{equation}
We use the terminology $L^{+}$-orbit of $\rF$ for an orbit that contains $L^{+}$ and likewise an $L^{-}$ -orbit. 
An orbit that contains neither $L^{+}$ or $L^{-}$ will be called a non-$L^{\pm}$ orbit. Because of the scale invariance of $\rF$, it suffices to study the orbit of $(0,1)$ to find an $L^{+}$ orbit and the orbit of $(0,-1)$ to find an $L^{-}$ orbit. As a result, we will often identify $L^{+}$, respectively $L^{-}$,
with $(0,1)$, respectively $(0,-1)$. 
More generally, the scale invariance of $\rF$ allows us to talk interchangeably about an orbit of points in the plane and the associated orbit of rays where each point is embedded into its position vector from the origin.

An orbit of $\rF$ is typically made from patching together \textbf{orbit segments} that occupy the domain of 
$\rM_a =\left(\begin{smallmatrix}a &-1\\ 1&0\end{smallmatrix}\right)$ in the right half plane with orbit segments that switch to occupy the domain of $\rM_b$ in the left half plane, before repeating this alternating behaviour. For reference, we denote these domains:
\begin{equation} \label{eq:domA}
\bd(\rM_a):=\bigl\{(x,y)\in\R^2\,:\,-\frac{\pi}{2}\leqslant 
 \arctan_2(y,x)<\frac{\pi}{2}\bigr\}, 
\quad \bd(\rM_b):= \R^2 \setminus \bd(\rM_a).
\end{equation}
In proposition \ref{prop:legalA} below, we study the dynamics within a single orbit segment in the right half plane (for ray counting in that result, 
note that $L^{-} \in \bd(\rM_a)$ but $L^{+} \in \bd(\rM_b))$.

Define the eigenvalues of $\rM_a$ by
\begin{equation} \label{eq:eigenv}
\lambda=\lambda_{\pm}= \frac{a}{2} \pm \sqrt{ \biggl(\frac{a}{2}\biggr)^2 - 1}, \quad \lambda_{+}\,\lambda_{-}=1.
\end{equation}
When $|a| \leqslant 2$, the elliptic case, we shall make use of the quantity
\begin{equation}\label{eq:Segments1}
 \zeta_j:=2\cos(\pi/j),\quad  j=1,2,\ldots.
\end{equation}
When $|a| > 2$, the hyperbolic case, the real eigenvectors associated to $\lambda_{\pm}$ are denoted $V_{\pm}$ and have slope $\lambda_{\mp}$.
From \eqref{eq:M^n} in the appendix, the powers of $\rM_a$ can be expressed in terms of polynomials $\cp_{n}(a)$ in $a$ of degree $n-1$, with $\cp_{n}(a)={\bar U}_{n-1}(a/2)$, where $\bar U_n$ is the $n$th Chebyshev polynomial of the second kind. The polynomials $U_n$ satisfy
\begin{equation}\label{eq:Chebysheva}
\cp_{-1}(a)=-1, \quad \cp_{0}(a)=0, \hskip 30pt 
\cp_{n+1}(a)=a\, \cp_{n}(a) -\cp_{n-1}(a), \quad n\geqslant 0.
\end{equation}
From \eqref{eq:M^n}, the forward images of $(1,0)=\rM_a(0,-1)$ by $\rM_
a^{j}$ for $j \geqslant 0$, as long as they remain in $\bd(\rM_a)$, are the rays
\begin{equation} \label{eq:ManRay}
(U_{j+1}(a), U_{j}(a)) \in \bd(\rM_a), \quad \mbox{with slope}\quad m_j:=U_{j}(a)/U_{j+1}(a).  
\end{equation}
Likewise, by reversibility, the images of $\rR\,(1,0)=(0,1)$ by $\rM_
a^{-j}$ are the rays
\begin{equation}  \label{eq:SecRay}
{\bf r}_j(a):= (U_{j}(a), U_{j+1}(a)) \in \bd(\rM_a),  \quad \mbox{with slope } 1/m_j.  
\end{equation}
The Mobius transformation that relates the slope $m$ of an initial ray in $\bd(\rM_a)$ to the slope $\mu_a(m)$ of its image is:
\begin{equation}\label{eq:Moba}
  m \mapsto \mu_a(m)=\frac{1}{a-m},  \quad m \in \R \cup \infty. 
\end{equation}
It inherits a reversing symmetry from $\rR$, given by
\begin{equation} \label{eq:defrev}
m \mapsto \frac{1}{m}.
\end{equation}

We also note that order preservation of the induced circle map corresponding 
to $\rM_a$ (or $\rM_b$) means that if one ray is obtained from another by an anti-clockwise rotation then the same is true of their images.  We use the anti-clockwise direction to order rays in the orbit segment from `first' to `last'.  

\begin{proposition} \label{prop:legalA}
Consider the dynamics of $\rM_a =\left(\begin{smallmatrix}a &-1\\1&0\end{smallmatrix}\right)$ in $\bd(\rM_a)$ of \eqref{eq:domA}.
\begin{enumerate}
\vspace*{-10pt}
\item [i)] For all $a > 2$, all rays between $L^{-}$ and the contracting 
eigenvector ray $V_{-}\,\cap\, \bd(\rM_a)$ with slope $\lambda_{+} > 1$ 
converge onto the expanding eigenvector ray $V_{+}\, \cap\, \bd(\rM_a)$ 
with slope $\lambda_{-} < 1$ ---the points on these rays escape to 
infinity along this direction; the rays between $V_{-}\, \cap\, \bd(\rM_a)$  and $L^{+}$ eventually escape to  $\bd(\rM_b)$.
\item [ii)] For any $\kappa\geqslant 2$, we have $a=\zeta_\kappa$ of 
\eqref{eq:Segments1} if and only if there is a segment of a symmetric orbit 
that goes from $L^{-}$ to $L^{+}$ with $\kappa$ rays in $\bd(\rM_a)$ and 
${\rM_a}^{\kappa}= -Id$.
\item [iii)] When $a \in (\zeta_{\kappa-1},\zeta_{\kappa}),\kappa\geqslant 2$,
the possible orbit segments of $\rM_a$ in $\bd(\rM_a)$ comprise a sequence of 
rays rotating anti-clockwise from the fourth quadrant into the first quadrant 
and eventually escaping to the left half-plane.  
The possibilities are:
\begin{enumerate}
\item  [a)] a $L^{+}$-orbit segment with $\kappa-1$ rays and first ray 
$(U_{\kappa-1}(a), U_{\kappa}(a))$ in the interior of the fourth quadrant;
\item [b)]  a $L^{-}$-orbit segment with $\kappa$ rays comprising 
$L^{-}=(0,-1)$ and $(1,0)$ and the $\rR$-image of the $L^{+}=(0,1)$-orbit 
segment contained within the first quadrant.
\item [c)] a non-$L^{\pm}$ orbit segment with $\kappa-1$ rays if the 
first ray has slope $m$ with $U_{\kappa}(a) / U_{\kappa-1}(a) < m < 0$ 
and $\kappa$ rays if the first ray has slope 
$m$ with $m<U_{\kappa}(a)/U_{\kappa-1}(a)<0$.
\end{enumerate}
\item [iv)] The rays ${\bf r}_j(a)$ of \eqref{eq:SecRay} that 
delineate the sectors in item iii a) rotate anti-clockwise as $a$ increases.
\end{enumerate}
\end{proposition}
\proof 
Firstly, we make a general comment about orbit segments. 
Since $(0,-1)$ maps anti-clockwise to $(1,0)$, independent of $a$, any orbit segment in $\bd(\rM_a)$ not starting parallel to $(0,-1)$ must have a single ray, its first, in the interior of the fourth quadrant (which is the image under $\rM_b$ of the last ray in the preceding orbit segment).  
This first ray with positive $x$-coordinate has an image with positive $y$-coordinate since $y'=x$ in $\rF$. 
On the other hand, the image of $(0,1)$ by $\rM_a^{-1}$ is $(1,a)$ and the 
last ray of any orbit segment of $\rM_a$ that escapes to the left half plane 
must be the single ray of the orbit segment located in the wedge sector 
determined by $(1,a)$ and $L^{+}$. 

We prove i). When $a > 2$, $U_{j+1}(a) > U_{j}(a) \geqslant j$ 
for $j \geqslant 1$ 
so \eqref{eq:ManRay} and order preservation implies the forward orbit 
of $(1,0)$, and hence that of any initial ray in the fourth quadrant, is confined thereafter to the first quadrant. 
Analysis of \eqref{eq:Moba} confirms that the fixed point at $m_{-}=\lambda_{-}$ of \eqref{eq:eigenv} is attracting on $[-\infty,\lambda_{+})$. 
The fixed point at $m_{+}=\lambda_{+}$ is repelling and the interval $(\lambda_{+}, \infty)$ eventually is mapped to $m < 0$ (which in this case corresponds to rays moving to the second quadrant). 
Of course, the fixed points 
$m_{\pm}=\lambda_{\pm}$ correspond to the eigenvectors of $\rM_a$.  

We prove ii). For $\rM_a^{-j}\,(0,1)=(U_{j}(a), U_{j+1}(a))$, requiring 
$U_j(a) > 0$ for $1 \leqslant j < \kappa$ and $U_{\kappa}(a) = 0$ produces 
$\kappa-1$ rays clockwise in the first quadrant, the last one parallel 
to $(1,0)$. 
Since the pre-image of $(1,0)$ is $(0,-1)$, we obtain $\kappa$ rays in total. 
For $\kappa=2$, the condition on $a$ for this orbit is that $U_2(a)=a=0$.  
For $\kappa > 2$, the condition $U_j(a)>0$ for $1 \leqslant j<\kappa$ and 
$U_{\kappa}(a) = 0$ finds the right-most root of $U_{\kappa}$, namely 
$a=\zeta_\kappa$. 
This follows from the well known roots of the Chebyshev polynomials $\bar U_n$,
as does the fact that $U_{\kappa-1}(\zeta_\kappa)=1$. 
Thus the $(\kappa-1)$th ray is exactly $(1,0)$, not just parallel to it.  
Substituting $U_{\kappa-1}(\zeta_\kappa)=1$, $U_{\kappa}(\zeta_\kappa) = 0$ and 
$U_{\kappa+1}(\zeta_\kappa)=-1$ into \eqref{eq:M^n} gives ${\rM_a}^{\kappa}= -Id$.  
Since $(1,0)=\rF(0,-1)$ and $(0,1)$ are mapped to one another by $\rR$, the finite orbit segment from $L^-$ to $L^+$ is part of a symmetric orbit.

iii) For an $L^{+}$-orbit segment of $\kappa-1$ rays, we require from 
\eqref{eq:SecRay} that  $U_j(a) > 0$ for $1 \leqslant j <  \kappa$ and $U_{\kappa}(a) < 0$, whence $a \in (\zeta_{\kappa-1},\zeta_{\kappa})$, giving part a).  
These $\kappa-1$ rays, together with $L^{+}$, delineate $\kappa$ sectors in $\bd(\rM_a)$.  
Taking any ray internal to the first sector bounded by $L^{-}$ and $(U_{\kappa-1}(a), U_{\kappa}(a))$ and iterating gives an orbit segment of $\kappa$ rays in $\bd(\rM_a)$, with one ray in each sector by order preservation. The first ray has slope $m:  m <  U_{\kappa}(a) / U_{\kappa-1}(a) < 0$.  
Taking any ray internal to the second sector starting at $(U_{\kappa-1}(a), U_{\kappa}(a))$ that is also in the fourth quadrant, with slope $m:  U_{\kappa}(a) / U_{\kappa-1}(a) < m < 0$, gives an orbit segment of $\kappa-1$ rays.  
This proves part c). 
A $L^{-}$ orbit segment of $\kappa$ rays begins with $L^{-}$ and then a ray parallel to $(1,0)=\rR\,(0,1)$.  From \eqref{eq:ManRay} and \eqref{eq:SecRay}, the forward orbit of $(1,0)$ contained in the first quadrant is the $\rR$-image of the backwards orbit of $(0,1)$ contained in the first quadrant, which from b) comprises $\kappa-1$ rays including now $(0, 1)$ itself.   

We prove iv). 
The inverse of \eqref{eq:Moba} is
$m \mapsto {\mu_a}^{-1}(m)=a - \frac{1}{m}$ for $m \in \R \cup \infty$, and it gives the slope $m(\rF^{-1} z)=\mu_a^{-1}(m)$ in terms of the slope $m=m(z)$ of a ray $z \in \R^2$.  We apply this to the rays ${\bf r}_j(a)$ of \eqref{eq:SecRay} that delineate the sectors in proposition \ref{prop:legalA} iii) a).  
Defining $r_j(a):=m({\bf r}_j(a))$, we have
$$r_{j+1}(a)=a -  \frac{1}{r_j(a)} \implies  r_{j+1}^{'}= 1 + \frac{r_j^{'}}{{r_j}^2}, $$
where the prime denotes differentiation with respect to $a$.   
Since $r_0=\infty, r^{'}_0=0$ and $r_1=a, r_1^{'}=1$, then by induction,
we see $r_{j}^{'} \geqslant 1$ for $j \geqslant 1$. Hence as $a$ increases,  ${\bf r}_j(a)$ moves anti-clockwise.
\endproof

There is the obvious version of the above for $\rM_b$.

This first dynamical analysis allows some preliminary bounds on rotation number and density in parameter space.
\begin{proposition} \label{prop:easybounds}
Let $\kappa \geqslant 2$. If $\a \in (\zeta_{\kappa-1},\zeta_{\kappa}]$, 
then the number of rays in each orbit segment in $\bd(\rM_a)$ is
$\kappa-1$ or $\kappa$. 
If $\b \in (\zeta_{\ell-1}, \zeta_{\ell}]$, $\ell \geqslant 2$, the analogous 
result holds in $\bd(\rM_b)$ and in addition:
$$
(i)\qad \theta(\a,\b) \in \bigl[\frac{1}{\kappa+\ell},\frac{1}{\kappa+\ell-2}\bigr]; 
\qquad
(ii)\qad \rho(\a,\b) \in  \bigl[\frac{\kappa-1}{\kappa+\ell-1},\frac{\kappa}{\kappa+\ell-1}\bigr].
$$
\end{proposition}
\proof
The statement on the number of rays versus parameter value follows directly from the list of possibilities in proposition \ref{prop:legalA}. 
If we use $N_a^i=N_a^i(\a,\b)$, $i \in \Z$, to label the bi-infinite sequence 
of the number of rays in $\bd(M_a)$ indexed by revolution $i$, and similarly 
for $N_b^i$, then
\begin{eqnarray*}
\theta(\a,\b) &=& \lim_{k \to \infty}\frac{2k+1}{\sum_{i=-k}^{i=k} N_a^i+N_b^i}, \\
\rho(\a,\b) &=&\lim_{k \to \infty}\frac{\sum_{i=-k}^{i=k} N_a^i}{\sum_{i=-k}^{i=k} N_a^i + N_b^i}\\
&=& \lim_{k \to \infty}  \bigl(1+ \frac{\sum_{i=-k}^{i=k} N_b^i}{\sum_{i=-k}^{i=k} N_a^i}\bigr)^{-1}.
\end{eqnarray*}
Using the bounds $(2k+1)\,(\kappa-1)  \leqslant  \sum_{i=-k}^{i=k} N_a^i 
 \leqslant (2k+1)\, \kappa$ etc., gives the results. 
Also we observe
$$ 
\rho(\a,\b) = K(\a)\,\theta(\a,\b),
$$
where $K(\a):=\lim_{k \to \infty} \frac{\sum_{i=-k}^{i=k} N_a^i}{2k+1}$ 
and  $\kappa-1 \leqslant K \leqslant \kappa$.
\endproof

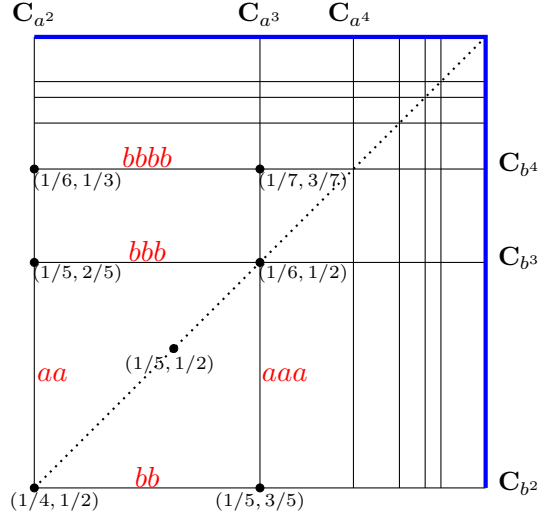
\begin{figure}[tb]
\centering
\begin{tikzpicture}[scale=0.3]
\draw[black, thin] (0,0) -- (20,0);
\draw[black, thin] (0,0) -- (0,20);
\draw[black, thin] (0,10) -- (20,10);
\draw[black, thin] (10,0) -- (10,20);
\draw[black, thin] (14.14213562,0) -- (14.14213562,20); 
\draw[black, thin] (0,14.14213562) -- (20,14.14213562,0); 
\draw[black, thin] (16.18033989,0) -- (16.18033989,20);
\draw[black, thin] (0,16.18033989) -- (20,16.18033989);
\draw[black, thin] (17.32050808,0) -- (17.32050808,20);
\draw[black, thin] (0,17.32050808) -- (20,17.32050808);
\draw[black, thin] (18.01937736,0) -- (18.01937736,20);
\draw[black, thin] (0,18.01937736) --(20,18.01937736);
\draw[black, dotted, thick] (0,0) -- (20,20);
\draw[blue, ultra thick] (0,20) -- (20,20) -- (20,0);


\node at (21.5,0.2) {\small $\bC_{b^2}$};
\node at (21.5,10.2) {\small $\bC_{b^3}$};
\node at (21.5,14.34) {\small $\bC_{b^4}$};

\node at (0,21.0) {\small $\bC_{a^2}$};
\node at (10,21.0) {\small $\bC_{a^3}$};
\node at (14,21.0) {\small $\bC_{a^4}$};

\node[red] at (5,0.5) {$bb$};
\node[red] at (5,10.5) {$bbb$};
\node[red] at (5,14.7) {$bbbb$};
\node[red] at (0.8,5) {$aa$};
\node[red] at (11.1,5) {$aaa$};

\filldraw[black] (0,0) circle (5pt);
\node[black] at (0.9,-0.6) {\tiny $(1/4,1/2)$};
\filldraw[black] (10,0) circle (5pt);
\node[black] at (10,-0.6) {\tiny $(1/5,3/5)$};
\filldraw[black] (0,10) circle (5pt);
\node[black] at (1.9,9.5) {\tiny $(1/5,2/5)$};
\filldraw[black] (6.180339876,6.180339876) circle (5pt);
\node[black] at (6,5.5) {\tiny $(1/5,1/2)$};
\filldraw[black] (10,10) circle (5pt);
\node[black] at (11.9,9.5) {\tiny $(1/6,1/2)$};
\filldraw[black] (0,14.14213562) circle (5pt);
\node[black] at (1.9,13.6) {\tiny $(1/6,1/3)$};
\filldraw[black] (10,14.14213562) circle (5pt);
\node[black] at (11.9,13.6) {\tiny $(1/7,3/7)$};

\end{tikzpicture}
\caption{\label{fig:RotDomains}\rm\small
Partition of the first quadrant $[0,2]^2$ of the parameter space $\cA$ 
into rotational domains, determined by the critical curves 
of proposition \ref{prop:Segments}.
The boundary of $\Theta(0)$ is represented in blue, and 
positive and negative critical words in red. 
The rational pairs near highlighted points represent the
value of rotational parameters at those points.
}
\end{figure}

\section{Critical curves}\label{section:CriticalCurves}

A \textbf{boundary parameter} is a pair $(\a,\b)\in\cA$ for which 
there is an $\rF$-orbit (a \textbf{boundary orbit}) that contains two (not necessarily distinct)
boundary rays, namely the positive and negative ordinate semi-axes
$L^+$ and $L^-$ of \eqref{eq:defLpm}.
Such a dynamical condition is accompanied by a finite \textbf{boundary word} 
$w$, which encodes the itinerary between two boundary rays; this part of
the boundary orbit will be called a \textbf{boundary segment}.
Because the two components of the map $\rF$ coincide on boundary rays, 
changing any letter that corresponds to a boundary ray has no effect
on the orbit. (This always applies to the first letter of a 
boundary word.) The resulting ambiguity defines an equivalence relation 
on words, whereby we write $w\sim w'$ to indicate that $w'$ is
obtained from $w$ by any change in the symbols allocated to a boundary ray.
A boundary word which is equivalent, but not equal, to the symbolic dynamics 
of the map $\rF$ is said to be \textbf{improper}. These words 
appear at the intersections of curves (section \ref{section:Intersections}).

The \textbf{rank} $\iota$ of a boundary word $w$ is given by 
$\iota=1+|w|_{ab}+|w|_{ba}$  [cf.~definition (\ref{eq:ThetaRho})].
If $w$ encodes the symbolic dynamics of $\rF$, then $\iota$ is the number of 
half-loops that occur in the transit between the initial and final 
boundary rays. 
A boundary word begins and ends with the same symbol
in the odd-rank cases 
($L^+\mapsto L^-$ or $L^-\mapsto L^+$)
and with distinct symbols in the even-rank case
($L^+\mapsto L^+$ or $L^-\mapsto L^-$).
The \textbf{sign} of a boundary word is positive if the word begins with 
the symbol $a$ (the orbit starts from $L^-$), and negative otherwise.

If the rank of $w$ is even, then the critical curve is the parametric locus of a 
periodic orbit having (at least) one point on the partition boundary. 
These are the periodic $\beta$-orbits in \cite[section 3.1]{Garcia-MoratoEtAl},
in which case the length and rank of the boundary word are, respectively, the 
denominator and twice the numerator of the rotation number.
If an even rank curve supports both positive and negative boundary segments,
then we call it an \textbf{axis} (of a resonance). 

If the rank of a word $w$ is odd, then the set of parameters for which 
$w$ is a boundary word will be called a \textbf{critical curve}, 
examples of which are given in \cite[examples 4.1--3]{LagariasRainsII} 
(see proposition \ref{prop:Segments} below).  The associated boundary orbit of $\rF$ will be called a 
\textbf{critical orbit} in this case.  

We see that the orbits of proposition~\ref{prop:legalA} ii) are critical orbits of rank 1 with critical curve $\a=\zeta_\kappa$ and boundary word $a^\kappa$ and this is true independently of the value of $\b$ and the dynamics in 
$\bd(\rM_b)$. This is also the family i) mentioned in the introduction, which will be used
in section \ref{section:RotationalDomains}.  
The following result extends example 3.2 in \cite{LagariasRainsI} 
(cf.~also propositions \ref{prop:legalA} and \ref{prop:easybounds} 
of the previous section):

\begin{proposition} \label{prop:Segments}
For $\kappa=2,3,\ldots$, the set
\begin{equation}\label{eq:Segments}
\bC_{a^\kappa}:=\{(\zeta_\kappa,\b) \,:\, \b\in\R\},
\end{equation}
with $\zeta_\kappa$ of \eqref{eq:Segments1} has rotational parameters
\begin{equation} \label{eq:ThetaRhoSegments}
\rho(\zeta_\kappa,\b)=\kappa\theta(\zeta_\kappa,\b),
\hskip 40pt
\theta(\zeta_\kappa,\b)=\begin{cases}
0&\b>2\\
\noalign{\vskip 3pt}
\displaystyle \frac{\arccos(\b/2)}{\kappa\arccos(\b/2)+\pi}& |\b|\leqslant 2\\
\noalign{\vskip 3pt}
\displaystyle\frac{1}{\kappa+1}& \b<-2.
\end{cases}
\end{equation}
Exchanging $\a$ and $\b$ yields the twin sets
$\bC_{b^\ell}:=\{(\a,\zeta_\ell)\,:\, \a \in\R\}$,
with the corresponding formula for the rotation number and
$\rho(\a,\zeta_\ell)=1-\rho(\zeta_\ell,\a)$.
\end{proposition}

\proof 
The rotation number in (\ref{eq:ThetaRhoSegments}) for $|\b|\leqslant 2$
follows from example 3.2 and theorem 2.1 in \cite{LagariasRainsI}.  
The value of $\theta$ for $\b>2$ follows from (the $\b$-part of)
proposition \ref{prop:legalA} and the proof of proposition 
\ref{prop:easybounds}.
If $\b \leqslant -2$, then by (\ref{eq:ThetaRhoSegments}) we have
$\theta(\zeta_\kappa,-2)={1}/{(\kappa+1)}$, and by
\cite[theorem 2.1]{LagariasRainsI}, $\theta(\a,\b)$ is nonincreasing in $\b$,
so that $\theta(\zeta_\kappa,\b) \geqslant 1/(\kappa+1)$.
However, (ii) of the same theorem implies $0 \leqslant \theta(\zeta_\kappa,\b)
 \leqslant 1/(\kappa+1)$ for any $\b$, and thus we have 
$\theta(\zeta_\kappa,\b)=1/(\kappa+1)$ for $\b \leqslant -2$.

The density in (\ref{eq:ThetaRhoSegments}) follows from $K(\a)$ defined 
in the proof of proposition \ref{prop:easybounds}, noting from 
proposition \ref{prop:legalA} ii) that for 
$\a=\zeta_\kappa$, we have $N_\a^i=\kappa$ for all $i$, 
whence $K(\zeta_\kappa)=\kappa$.
Exchanging $\a$ and $\b$ yields the same analysis; the obvious modification
for the density results from the congruence $\rho\equiv \kappa\theta\mod{1}$
with negative $\kappa$.
\endproof

Note that in the interval $|\b|\leqslant 2$ ($|\a|\leqslant 2$) these 
components do not intersect the interior of any resonance, by virtue
of the fact that the rotation number is nowhere constant.
The points $(\zeta_\kappa,\pm2)$ of $\bC_{a^\kappa}$ lie
on the boundary of the resonances with rotation number $0$ 
and $1/(\kappa+1)$.

The critical curves are our main object of study, and the closure of
the set of boundary parameters which belong to some critical curve
will be called the \textit{critical set} $\Xi\subset \cA$.
The simplest components of the critical set are those described in
proposition \ref{prop:Segments}.
Indeed, for $\kappa\geqslant 2$, the positive word 
$a^\kappa=\underbrace{a\cdots a}_\kappa$ is critical of rank one,
with components $\bC_{a^\kappa}$.
Likewise, $b^\ell$ is negative of rank one, with component $\bC_{b^\ell}$.

For later use, we collect the following which are known or easily verified:
\begin{lemma} \label{lemma:easy}
For the map $\rF=\rF_{ab}$ of \eqref{eq:LR} and its reversing symmetry $\rR$ of \eqref{eq:defR}, we have:
\begin{enumerate}
\item [i)] From reversibility, the $\rR$-image of the forward (backward) 
orbit of $(0,y)$ is the backward (forward) orbit of 
$\rR\,(0,y)=(y,0)=\rF(0,-y)$, i.e.,
$$ 
\rF^{-k+1}\,(0, -y)=\rR\,\rF^{k}\,(0,y),\quad k \in \Z,\;y\in\R^{+}.
$$
\item  [ii)] A critical orbit of odd rank exists if and only if $(0,1)$ and $\rR\,(0,1)=(1,0)=\rF(0,-1)$ are in the same orbit.
Since the point $(0,1)$ and its image under $R$ are in the same orbit,  the critical orbit is symmetric.
\item [iii)] With $T: (x,y) \mapsto (-x,-y)$, we have
$$   \rF_{ab} = T \circ \rF_{ba} \circ T^{-1}. $$
In particular, it suffices to study critical orbits where $(0,-1)$ goes to $(0,1)$ in forward time as those from $(0,1)$ to $(0,-1)$ for $F_{ab}$ map to the former for $F_{ba}$.
\item [iv)] There can be no critical orbit (of rank $>1)$ from $(0,-1)$ to $(0,1)$ if $a > 2$ ($b > 2$).
\item [v)] We have  $\rF (0,\pm1)=(\mp1,0)$ together with
$\rF (1,0)=(a,1),\; \rF (-1,0)=(-b,-1)$ and  $\rF^{-1} (0,1)=(1,a), \rF^{-1} (0,-1)=(-1,-b)$. Additionally:
 \newline (a) each orbit segment in $\bd(\rM_a)$ $(\bd(\rM_b))$ only has at most one point in the interior of the fourth (second) quadrant, with its image in the first or second (third or fourth) quadrant;
\newline (b) if $a,b \in (-2,2)$, and if an orbit segment has a point $P$ in the interior of the first (third) quadrant, then some image $\rF^{j}\, P$ with $j \ge 1$ is in the second (fourth) quadrant and some pre-image $\rF^{-k}\, P$ with $k \ge 1$ is in the fourth (second) quadrant.  
\end{enumerate}
\end{lemma}
\proof
Property i) is a standard result for reversible maps \cite{LambRoberts}.  
Property ii) is from \cite[theorem 3.4]{LagariasRainsI}. 
Property iii) is from \cite[equation (2.4)]{LagariasRainsI}. 
Property iv) follows from Proposition \ref{prop:legalA} i) and its 
corresponding version for $\rM_b$.
Property v) follows from $\rF$ being order preserving and $\rF: (x,y) \mapsto (x',y')$ having $y'=x$ and $\rF^{-1}: (x',y') \mapsto (x,y)$ having $x=y'$.
\endproof

\subsection{Algebraic properties}\label{section:AlgebraicProperties}

Requiring that an orbit contain two boundary rays leads to algebraic 
curves over $\mathbb{Q}$. We introduce the necessary formalism, and 
establish some properties of these curves.

Let
\begin{equation}\label{eq:M}
\rM_x =\left(\begin{matrix}x &-1\\ 1&0\end{matrix}\right)
\end{equation}
where $x$ is an indeterminate and consider the 
finite word  
$$
w=w_0w_1\cdots w_{t-1},$$
where we regard, for now, each letter of the word $w$ as an indeterminate 
(in the case of $\rF$, we will later specialise to $w_i \in \{a,b\}$).
We form a product of matrices of type (\ref{eq:M}) as follows:
\begin{equation}\label{eq:Mw}
\rM_{\varnothing}=\mathrm{Id},\qquad
\rM_{[w]}=\rM_{w_{t-1}}\rM_{w_{t-2}}\cdots \rM_{w_0},
\end{equation}
where $\varnothing$ denotes the empty word.

The entries of $\rM_{[w]}$ are polynomials in the indeterminates 
$w_0,w_1\cdots, w_{t-1}$.
Define the recursive sequence of polynomials
$Q_t=Q_t(w_1,\ldots,w_{t-1}) \in \Z[w_0,w_1,\ldots,w_{t-1}]$ 
by the three-term recurrence:
\begin{equation}\label{eq:PQ}
\begin{array}{llll}
Q_{-1}=-1 & \quad Q_0=0 & \hskip 30pt 
      Q_{t+1}=w_{t}Q_t-Q_{t-1}& \quad t\geqslant 0.
\end{array}
\end{equation}
Hence 
\begin{eqnarray*}
Q_1&=&1,
\\
Q_2(w_1) &= &w_1,
\\
Q_3(w_1,w_2)&=&w_1w_2-1,
\\
Q_4(w_1,w_2,w_3) & = & w_1w_2w_3 - w_1 - w_3,
\\
Q_5(w_1,w_2,w_3,w_4) &= & w_1w_2w_3w_4- w_1 w_4 - w_3 w_4 - w_1 w_2+1.
\end{eqnarray*}
Alternatively, we can use the tridiagonal determinant representation:
\begin{equation} \label{eq:continuant}
Q_{t}(w_1,\ldots,w_{t-1}):= 
\begin{vmatrix}
w_1 & 1 & 0 & \cdots & 0 \\
1 & w_2 & 1 & & \vdots \\
0 & \ddots & \ddots & \ddots & 0 \\
\vdots &  & 1 & w_{t-2} & 1 \\
0 & \cdots & 0 & 1 & w_{t-1} \\
\end{vmatrix}, \quad t \ge2,
\end{equation} 
with the recurrence relation \eqref{eq:PQ} following by expanding the 
determinant along the last row.  
Note we can also run the recurrence \eqref{eq:PQ} backwards using 
indeterminates $w_{-i}$ to find:
\begin{equation} \label{eq:NegInd}
Q_{-t} = -Q_t(w_i \mapsto w_{-i}) \qquad t \geqslant 0.
\end{equation}

With $Q_t$ defined, it follows by induction:
\begin{equation}\label{eq:MwReduced}
\rM_{[w_1\cdots w_{t-1}]} =
\begin{pmatrix}
Q_{t}(w_1,\ldots,w_{t-1}) &-Q_{t-1}(w_2,\ldots,w_{t-1}) \\
Q_{t-1}(w_1,\ldots,w_{t-2}) & -Q_{t-2}(w_2,\ldots,w_{t-2})
\end{pmatrix},
\end{equation}
and 
\begin{equation}\label{eq:Mws}
\rM_{[w_0\cdots w_{t-1}]} =
\rM_{[w_1\cdots w_{t-1}]}\rM_{[w_0]} =
\begin{pmatrix}
  w_0Q_{t}-Q_{t-1} & -Q_{t}\\
  w_0Q_{t-1}-Q_{t-2} & -Q_{t-1}
\end{pmatrix}\qquad t\geqslant 1.
\end{equation}

The above results follow from the theory of {\em continuants} 
(or {\em continuant polynomials}), a name given to the tridiagonal 
determinant \eqref{eq:continuant}. 
Continuants and their properties were studied by Euler 
\cite[Section 6.7]{GrahamEtAl} in connection with generalised 
continued fractions involving arbitrary real (or complex) numbers 
instead of integers, see proposition \ref{prop:ContinuedFractions} i) below. 
The  associated three-term recurrence for Euler's original continuant 
is for the polynomial sequence $q_t$ satisfying $q_{-1}=-1$, $q_0=0$ 
and $q_{t+1}=w_{t}q_t+q_{t-1}$.  
The polynomial $Q_t$ of \eqref{eq:PQ} belongs to a class of generalisations 
of the classic continuant, variously called  signed continuant polynomials 
or generalised Chebyshev polynomials,  which have arisen in the study of 
cluster algebras and frieze patterns \cite{MGO, B-MEtAl}. 
The relationship between $Q_t$ and $q_t$ is 
$Q_t(w_1,\ldots,w_{t-1}) = (-i)^{t-1}\,q_t(i\,w_1,\ldots,i\,w_{t-1})$ 
with $i=\sqrt{-1}$.  

It follows from Euler that $Q_t$ can be generated as the sum of the 
product $w_1 w_2 \ldots w_{t-1}$, its leading term, together with 
all possible ways of writing the leading term again but striking 
out adjacent product pairs $w_i w_{i+1}$ and replacing such a  
pair with $-1$ (e.g., the four non-leading terms of $Q_5$ above are 
obtained by striking out, respectively, $w_2 w_3$, $w_1 w_2$, $w_3 w_4$  
and the two pairs $w_1 w_2$ and $w_3 w_4$). 
Note that from \eqref{eq:PQ} and \eqref{eq:Chebysheva},
\begin{equation} \label{eq:RedCheby}
Q_{n+1}(x,x,\ldots x) = \cp_{n+1}(x) = \bar U_n (x/2), 
\end{equation}
where $\bar U_n$ are the Chebyshev polynomials of the second kind ---see 
the Appendix.

We shall need the following properties of continuants, 
found in references \cite[Chapter 6.7]{GrahamEtAl}\cite[Section 5.1]{MGO}, 
\cite[Section 2.2]{B-MEtAl}:
\cite{GrahamEtAl,MGO,B-MEtAl}, 
\begin{proposition}[Continuant Polynomials] \label{prop:ContinuedFractions}
For all $t>1$, indeterminates $w_1,\ldots,w_{t-1} \in \R$ and integers $k,l 
\geqslant -1$, we have
\vspace*{-10pt}
\begin{enumerate}
\item [i)]
$\displaystyle
\frac{Q_t(w_1,\ldots,w_{t-1})}{Q_{t-1}(w_2,\ldots,w_{t-1})}=
                  \,w_1-\frac{1}{\displaystyle\strut
                  \,w_2-\frac{1}{\displaystyle\strut
                  \,\ddots -\frac{1}{w_{t-1}}}} 
$
\item [ii)]
$\displaystyle
Q_t(w_1,\ldots,w_{t-2},w_{t-1})=Q_t(w_{t-1},w_{t-2},\ldots,w_1)
$
\item [iii)]
$\displaystyle
Q_{k+l+1}(w_1,\ldots,w_{k+l})   =  \\
 Q_{k+1}(w_1,\ldots,w_{k})\, Q_{l+1}(w_{k+1},\ldots,w_{k+l}) - 
Q_{k}(w_1,\ldots,w_{k-1})\, Q_{l}(w_{k+2},\ldots,w_{k+l}).
$
\end{enumerate}
\end{proposition}
                     
We now specialise to the case of odd rank, for which 
the structure of words and curves is constrained by the 
\textit{reversibility} of the map.
The \textbf{reduced word} $\underline w$ of $w$ is defined as 
$\underline w=w_1\cdots w_{n-1}$ (again reversibility dictates considering to drop the first letter $w_0$ -- see below). We say that $\underline w$
is a \textit{palindrome} if
$$
w_i=w_{n-i}\qquad i=1,\ldots,n-1.
$$
For $\rF$, the structure of a (palindromic) reduced boundary word of odd rank can also be encoded by the (palindromic) integer exponent sequence $(i_k)$  of odd length $2m-1$, using $i_k \in \Z^+$, $1 \leqslant k \leqslant m$,  that we call its \textbf{block sequence}, e.g., 
\begin{equation} \label{eq:block}
 \underline w=a^{i_1}\, b^{i_2}\, a^{i_3}\, \ldots  b^{i_{m-1}}  a^{i_m}\,  b^{i_{m-1}}\,  \ldots a^{i_3}\, b^{i_2}\, a^{i_1}.
\end{equation}
In this form, the rank is now obvious, being $2m-1$, whereas 
\begin{equation} \label{eq:balance}
i_m+ 2\,\sum_{k=1}^{m-1} i_k = n-1. 
\end{equation}
We take the convention that $i_k$, $k$ odd, counts powers of $a$ if the associated word is positive as in \eqref{eq:block} (recall this means the word begins with $a$) and powers of $b$ if the word is negative.  We see for a positive palindromic word that the middle block comprises powers of $a$ when $m$ is odd and powers of $b$ when $m$ is even.  We recall from proposition \ref{prop:easybounds}  that necessarily:
\begin{equation} \label{eq:blockgeom}
\a \in (\zeta_{\kappa-1},\zeta_\kappa) \implies i_k \in \{\kappa-1,\kappa\}, \quad
\b \in (\zeta_{\ell-1},\zeta_\ell) \implies i_k \in \{\ell-1,\ell \}.
\end{equation}

The significance of palindromic words to boundary curves of the
map $\rF$ is established by the following result (it suffices to prove the case of a positive word noting ii) and iii) of lemma \ref{lemma:easy}).

\begin{theorem}\label{thm:PalindromeII}
Consider a non-periodic critical orbit of $\rF$ that contains $(0,-1)$ and $(0,1)$ and the associated positive boundary word $w$ of length $n$ that encodes the itinerary between them. We have: 
\vspace{-10pt}
\begin{enumerate}
\item [i)] The reduced boundary word of rank $2m-1$ is a palindrome \eqref{eq:block} of $2m-1$ blocks (equivalently its block sequence is a palindromic $(2m-1)$-tuple of positive integers)
\item [ii)] The integer $n$ is odd and $i_m$ is even if and only if the boundary segment intersects $Fix(\rR)$ once in $\bd(\rM_a)$  $(\bd(\rM_b))$, whence $m$ is odd (even). 
\item  [iii)] The integer $n$ is even and $i_m$ is odd if and only if the boundary segment intersects $Fix(\rF \circ \rR)$ once. The intersection is in $\bd(\rM_a)$  $(\bd(\rM_b))$ if $i_m=1$ and $m$ is even (odd) or if 
$i_m > 1$  and $m$ is odd (even).
\end{enumerate}
\end{theorem}
\vspace{-10pt}
\proof
i). From reversibility, for arbitrary $z \in \R^2$, we have $z'=\rF(z) \iff \rR z = \rF(\rR z')$.  We claim the symbol sequence of the $\rF$-orbit leaving the point $z\not=(0,\pm 1)$, developing to the right, is the same as the symbol sequence for the $\rF$-orbit arriving to $\rR z$, developing to the left.  
To see this, note the first symbol for the former is determined by $z$ and the first symbol for the latter is determined by $\rR z'$ and we claim that these symbols are the same. 
From proposition \ref{lemma:easy} v), we have that if $z$ is in the first or 
fourth quadrant, hence encoded with symbol $a$ for its forward image, 
then $z'$ is in the first or second quadrant, whence $\rR z'$ is in the first or fourth quadrant, and $\rR z'$ is encoded similarly to $z$. 
A similar result is true for $z$ in the second or third quadrant and the symbol $b$.  
The argument is then iterated, next with $z'$ and $\rR \rF(z')$.  If we take a critical orbit with $z=(1,0)=\rF(0,-1)$, we see the finite orbit connecting $z$ and $\rR z =(0,1)$ must have a palindromic symbol sequence, i.e., the reduced word is a palindrome as in \eqref{eq:block}.

We prove ii). We know a critical orbit is symmetric from lemma \ref{lemma:easy}.
In general \cite{LambRoberts}, for any $G$ that acts as a reversor of a reversible map $L$, so $L^{-1}=G^{-1} \circ L \circ G$,  an orbit is $G$-invariant  if and only if  $G\, z = L^{k_z} z$ with  $k_z \in \Z$ for each point $z \in \R^2$ of the orbit. The latter is true if and only if it is true for one point. Then either: 
\newline (a) $G\, L^{j_z}\,  z = L^{j_z}\, z$ if $k_z=2j_z$, i.e., $L^{j_z}\,z$ is fixed by $G$ ; or 
\newline (b) $LG\, L^{j_z+1} z = L^{j_z+1} z$  if $k_z=2\,j_z+1$, i.e.,  $L^{j_z+1} z$ is fixed by $LG$.
If a symmetric orbit is not periodic, it contains one point fixed by $G$ or one point fixed by $LG$.

We are interested in $L=\rF$ and $G=\rR$ and we can take the point $z=(1,0)$ 
and $k_z=n-1$. 
If $n$ is odd so $k_z$ is even, we have case (a). 
Hence the forward orbit of the midway point $\rF^{k_z/2}\, (1,0) \in Fix(\rR)$ is the reflection by $\rR$ of its backwards orbit and the slope of the rays in the forward and backwards orbit are reciprocals of each other from \eqref{eq:defrev}. 
The involution $\rR$ preserves the first and third quadrants and interchanges the second and the fourth. 
If $m > 1$ is even (odd), the middle block in \eqref{eq:block} is built from the letter $b$ (the letter $a$) and the number of points $i_m$ of this orbit segment in $\bd(M_b)$ ($\bd(M_a)$) must be even.  
It comprises $\rF^{k_z/2}\, (1,0)$ on the line $y=x$ and pairs of its forward and backward iterates in the interior of the third (first) quadrant, plus the necessary additional point in the second (fourth) quadrant guaranteed from lemma \ref{lemma:easy}.  
If $m=1$, we have only a single letter in the word and then $i_m$ in the reduced word is even (and the power is odd in the full word, cf.~proposition \ref{prop:legalA}).  
Hence $i_m$ is even as claimed, which can also be seen from \eqref{eq:balance}. 

We prove iii).
Now consider the case (b) above, i.e., $n$ is even so $k_z$ is odd and $q=\rF^{j_z+1}\, (1,0) \in Fix(\rF \rR)$. 
For $y >0$, the involution $\rF \rR: x'=-x + a y, y'=y$ is a horizontal reflection about the point $\frac{a}{2}\,y$; for $y <0$, it is a horizontal reflection about the point $\frac{b}{2}\,y$. 
We have $\rF^j \,q = \rF^j\, (\rF \rR\, q) = (\rF \rR) (\rF^{-j}\, q)$, so the forward and backward iterates of $q$ under $\rF$ at corresponding times in the upper half-plane and in the lower half-plane must be pairs under these reflections. 
As an example, suppose $q \in \bd(\rM_a)$ in the first quadrant. 
This means the single ray in the second quadrant when $\rF^r q$ first enters 
it (cf.~lemma \ref{lemma:easy}) must force the iterates of $q$ in the 
first quadrant to have $r-1$ rays to the left of $x=\frac{a}{2}\,y$ in forwards time, $r \geqslant 1$, and $r$ rays to the right in backward time.  
Counting $q$ itself and the single ray $\rF^{-(r+1)}\, q$ in the 
fourth quadrant gives $2r+1$ for the number $i_m$ of the letter $a$ 
in the middle block, whence $i_m$ is odd. 
Similar reasoning applies for the other possibilities.
\endproof

This result leads us to study continuant polynomials for palindromic words (interestingly, palindromic continuants were used in Smith's 1855 proof of the Fermat two-square theorem \cite{ClarkeEtAl}). For generality, we again let $w=w_0\cdots w_{n-1}$ be a word in $n$ letters (not just two letters $a$ and $b$). Let
\begin{equation}\label{eq:n2}
\n2=\lfloor n/2\rfloor.
\end{equation}
and build the particular polynomials in
$\Z[w_0,\ldots,w_{n-1}]$ from  (\ref{eq:PQ}):
\begin{equation}\label{eq:C}
C_w :=\begin{cases}
Q_{\n2+1}-Q_{\n2}& \mbox{$n$ odd}\\
Q_{\n2+1}-Q_{\n2-1}& \mbox{$n$ even}
\end{cases}
\hskip 35pt
\Cbar_w:=\begin{cases}
Q_{\n2+1}+Q_{\n2}& \mbox{$n$ odd}\\
Q_{\n2}& \mbox{$n$ even}.
\end{cases}
\end{equation}
For completeness, if the rank of $w$ is \textit{even}, we let 
$C_w=Q_n$ and $\Cbar_w=1$.

The following result collects some algebraic properties
of the $Q_t$ polynomials for palindromic words in any alphabet.

\begin{proposition} \label{prop:Palindrome} 
If $w_1\cdots w_{n-1}$ is a palindrome (in any alphabet), 
then the following holds for the polynomials $Q_t$ of \eqref{eq:PQ} with $t=0,\ldots,n$ (noting also \eqref{eq:NegInd} for (ii)):
\begin{eqnarray}
i)&&Q_n=C_w\Cbar_w
\label{eq:Factorisation}\\
\noalign{\vskip 5pt}
ii)&&Q_{n-t}-Q_t=C_w\times
\begin{cases}
Q_{t-\n2+1}(w_{\n2}=-1, w_{\n2+1}, \ldots, w_{t-1}) & \mbox{$n$ odd}\\
-Q_{t-\n2}(w_{\n2+1}, \ldots, w_{t-1}) & \mbox{$n$ even}.
\label{eq:QSymmetry}
\end{cases}
\end{eqnarray}
\end{proposition}

\proof i) From proposition \ref{prop:ContinuedFractions} ii) and  iii) 
with $l=t-1$, $k=n-t$ and the palindromic word, we have for $t=1,\ldots,n-1$:
\begin{equation}. \label{eq:Recursion}
\begin{split}
Q_{n}(w_1,\ldots,w_{n-1})   & =  \\
 Q_{n-t+1}(w_1,\ldots,w_{n-t})\, Q_{t}(w_{1},\ldots,w_{t-1})  &- 
Q_{n-t}(w_1,\ldots,w_{n-t-1})\, Q_{t-1}(w_{1},\ldots,w_{t-2}).
\end{split}
\end{equation}
For odd $n=2\n2+1$, we take  $t=\n2+1$, to obtain 
$Q_n=Q_{\n2+1}^2-Q_\n2^2=C_w\Cbar_w$. 
For even $n=2\n2$, we take $t=\n2=n-t$, giving 
$Q_n=Q_{\n2+1} Q_\n2 - Q_\n2 Q_{\n2-1}=C_w \Cbar_w$.

ii)  We prove the result by induction, in the first instance for $t \geqslant \n2$.  
Consider the polynomial $A_t=Q_{n-t}-Q_t$. If $n=2\n2+1$,
from (\ref{eq:C}) we find:
$$
A_{\n2}=Q_{\n2+1}-Q_\n2=C_w=C_w\, Q_1
\qquad
A_{\n2+1}=Q_{\n2}-Q_{\n2+1}=-C_w=C_w\, Q_2(w_{\n2}=-1).
$$
If $n=2\n2$, then
$$
A_{\n2}=Q_{\n2}-Q_{\n2}=0=C_w\, Q_0
\qquad
A_{\n2+1}=Q_{\n2-1}-Q_{\n2+1}=-C_w=C_w\, (-Q_1).
$$
The above data serves as the base case for induction.
Assume that (\ref{eq:QSymmetry}) holds for all $i$ in the
range $\n2\leqslant i\leqslant t$, for some $t  \geqslant \n2+1$.
Then
\begin{eqnarray*}
A_{t+1}
 &=&Q_{n-(t+1)}-Q_{t+1}=w_{n-t}Q_{n-t}-Q_{n-(t-1)}-w_tQ_t+Q_{t-1}\\
 &=&w_t\,(Q_{n-t}-Q_t)-(Q_{n-(t-1)}-Q_{t-1})= C_w\, (w_t S_{t} - S_{t-1})
 = C_w\, S_{t+1},
\end{eqnarray*}
where $S_t$ denotes the particular solutions in \eqref{eq:QSymmetry} for the respective cases of $n$ even and $n$ odd, which obviously satisfy $S_{t+1}=w_t S_{t} - S_{t-1}$.
This completes the induction for the range $t\geqslant\n2$.  To see that \eqref{eq:QSymmetry} holds also for $0 \leqslant t < \n2$, note that both sides become their negatives under $t \mapsto n-t$, using \eqref{eq:NegInd} on the right hand side.
\endproof

We now apply the above results to the dynamics of $\rF$. 
Let $w=w_0\cdots w_{n-1}$ be a boundary word -- necessarily $n \geqslant 2$ for a boundary word
since one matrix of the form (\ref{eq:M}) cannot map a boundary ray to a boundary ray.  Let
$x_{-1}=-\mathrm{sign}(w),x_0=0,\ldots,x_n$ be the corresponding orbit segment.
From the equation
$$
\begin{pmatrix} x_t\\ x_{t-1}\end{pmatrix}
=\rM_{[w_0\cdots w_{t-1}]}
\begin{pmatrix}0\\-\mathrm{sign}(w)\end{pmatrix}
$$
and \eqref{eq:Mws}, 
we obtain 
$$x_t=\mathrm{sign}(w)\,Q_{t}(w_1,\ldots,w_{t-1}), \quad t \geqslant 1.$$
An equation for the  boundary curve with word $w$ is obtained by
requiring that $x_n=0$:
\begin{equation}\label{eq:EqCurveAll}
Q_n(\a,\b)=0.
\end{equation}
This equation stores redundant information.
It expresses the fact that the image of one boundary ray under
the matrix product $\rM_{[w]}$ is another boundary ray, and
such an action may be realised without any reference to
the symbolic dynamics of the map $\rF$. 
As a result, the curve 
\begin{equation}\label{eq:CurveAll}
\cQ_n=\{(\a,\b)\in\R^2\,:\,Q_n(\a,\b)=0\}
\end{equation} 
has several branches, as we shall see below.
A parameter pair $c=(\a,\b)$ for which $w$ is equivalent\footnote{%
In the sense mentioned at the beginning of section \ref{section:CriticalCurves}.} 
to the symbolic word of a boundary segment of $\rF$ will be called a 
\textbf{legal point} of the curve, and a \textbf{legal branch} of 
the curve is one containing legal points.
If a point is not legal, then it may happen that the rank of
the word does not correspond to the number of half-turns performed
by the orbit segment. We call the latter the \textbf{orbital rank} of
the point $c$.

\goodbreak
\begin{theorem} \label{thm:Branches}
Let $\cQ_n$ be the curve (\ref{eq:CurveAll}), with reduced word $\underline w=w_1\cdots w_{n-1}$, $n \geqslant 2$.
Then
\vspace*{-10pt}
\begin{enumerate}
\item [i)] $\cQ_n$ has $n-1$ disjoint branches, 
of which precisely one is legal. On the legal branch,
rank and orbital rank coincide, and vice-versa.
\item [ii)] If the rank of $w$ is greater than 1, then each branch is 
represented by a decreasing function $b=b(a)$.
\item [iii)] $\cQ_n$ has $n-1$ asymptotes, $|\underline w|_b$ of which horizontal 
and $|\underline w|_a$ vertical, including multiplicities.
\end{enumerate}
\end{theorem}

\proof We prove i).
From (\ref{eq:RedCheby}), we see that
for $a=b$ we have $Q_n(a,a)=\bar U_{n-1}(a/2)$. 
Since the latter Chebyshev polynomial
has $n-1$ distinct real roots, the curve (\ref{eq:CurveAll}) 
intersects the line $a=b$ in $n-1$ distinct points, which are
\begin{equation}\label{eq:zeta}
\a=\b=\zeta_{n,j}=2\cos\left(\frac{\pi j}{n}\right),\qquad j=1,\ldots,n-1.
\end{equation}
[The number $\zeta_{n,1}$ above corresponds to 
$\zeta_n$ in (\ref{eq:Segments1}).]
Over that line, we have $\rM_{[w]}=\rM_a^{n}$, and 
as $\a$ decreases from $2$ to $-2$, the rotation number of $\rM_a$ 
increases monotonically from 0 to 1/2. At the value $\a=\zeta_{n,j}$
the image of a boundary ray will reach a boundary ray after $j$ 
half-turns, and therefore the orbit of $\rM_a$ will have the correct
rank precisely for $j=\iota$.

Now prolong each branch starting from the corresponding point (\ref{eq:zeta}).
Since the prolongation preserves the initial and final rays, 
as well as the orbital rank, distinct branches cannot intersect. 
Thus the boundary curve has a unique legal branch, namely that
where rank and orbital rank coincide, which the branch
containing the point $\a=\b=\zeta_{|w|,\iota(w)}$. 
The proof of i) is complete.

We prove ii). Let $c=(\a,\b)$ be a (finite) point on $\cQ_n$. 
We consider the 
$\rM_{[\underline w]}(c)$-orbit of the appropriate ray $L_1=(\pm1,0)$, 
the sign agreeing with that of $w$ (the notation $\rM_{[\underline w]}(c)$ refers to putting the specified parameter values $c=(\a,\b)$ into the matrix entries of \eqref{eq:MwReduced}).
We don't require $c$ to be legal, 
so the word $w$ may be unrelated to
the symbolic trajectory of the orbit of $\rF$.
Since the rank of $w$ is greater than 1, 
the orbit segment of $w$ will have at least one non-boundary 
ray acted upon by each matrix $\rM_a$ and $\rM_b$, where
$\rM_x$ is given in (\ref{eq:M}). 

One verifies that for any ray $L$ and parameters $\xi$ and $\epsilon$, 
the ray $\rM_{\xi+\epsilon} \,L$ is obtained from $\rM_{\xi}\,L$ by rotating 
clockwise if $\epsilon>0$, and anticlockwise if $\epsilon<0$, unless 
$L=L^\pm$, in which case the two rays 
coincide\footnote{This is theorem 3.2 (i) of \cite{LagariasRainsI}.}.
Considering that the circle map $\rf$ is orientation-preserving
\cite[theorem 3.1]{LagariasRainsI}, by repeating the above argument
it follows that for any $\epsilon>0$, both 
$\rM_{[\underline w]}(c+(\epsilon,0))\,L_1$ and
$\rM_{[\underline w]}(c+(0,\epsilon))\,L_1$ are obtained from 
the vertical boundary ray $\rM_{[\underline w]}(c)\,L_1$ via a clockwise rotation. 

From the above argument, it follows that the partial derivatives 
$\partial Q_n/\partial a$ and 
$\partial Q_n/\partial b$ are non-zero and agree in sign\footnote{They are
both positive if the rank is odd, and negative if the rank is even.}, 
that is, the tangent at any point of the curve has negative slope.
(This also shows that $\cQ_n$ has no isolated points.) 
We have proved ii).

We prove iii). Let $n=|w|$, $n_a=|\underline w|_a$, $n_b=|\underline w|_b$,
and 
\begin{equation}\label{eq:Cpolynomial}
Q_n(a,b)=a^{n_a}b^{n_b}
  + \sum_{k=0}^{n-2}\sum_{i=0}^k c_{i,k-i}a^ib^{k-i},
    \quad\mbox{for some}\quad c_{i,j}\in\Z.
\end{equation}
We rewrite the above as follows
\begin{equation}\label{eq:Asymptotes}
\begin{array}{rcll}
\displaystyle
\frac{1}{a^{n_a}}Q_n(a,b)&=&
\displaystyle
     b^{n_b}+\sum_{j=0}^{n_b-1}c_{n_a,j}b^{j}
    +O\bigl(\frac{1}{a}\bigr)&\quad |a|\to\infty\\
\displaystyle
\frac{1}{b^{n_b}}Q_n(a,b)&=&
\displaystyle
     a^{n_a}+\sum_{j=0}^{n_a-1}c_{j,n_b}a^{j}
    +O\bigl(\frac{1}{b}\bigr)&\quad |b|\to\infty.
\end{array}
\end{equation}
As $|a|$ or $|b|$ tends to infinity, the corresponding polynomial on 
the RHS of (\ref{eq:Asymptotes}) must vanish, each root giving the 
equation of an asymptote. This gives $n_a+n_b=n-1$ asymptotes, counting
multiplicities. Since the curve $\cQ_n$ has order $n-1$, by B\'ezout's theorem, 
it cannot have more than $n-1$ points on the line at infinity, so
there are no other asymptotes.
\endproof

Let us return to equation (\ref{eq:EqCurveAll}).
If the rank of $w$ is even, then $Q_{n}$ may be (and typically is)
irreducible, as in the case $w=a^2b^2$ for which $Q_4=ab^2-a-b$.
Thus the equation (\ref{eq:EqCurveAll}) is the minimal description
of a boundary curve of even rank, in general. 

For odd rank, we consider the factorisation 
(\ref{eq:Factorisation}), and replace (\ref{eq:CurveAll}) by
\begin{equation}\label{eq:CriticalCurve}
\cC_w=\{(a,b)\in\R^2\,:\,C_w(a,b)=0\},
\end{equation}
where $C_w$ is defined in (\ref{eq:C}).
This is justified as follows. From reversibility we have 
$x_t=x_{n-t}$ for all $t$. 
If $n$ is odd, then the point 
 $z_{\n2+1}=(x_{\n2+1},x_\n2)=\pm(Q_{\n2+1},Q_\n2)$ 
lies on the symmetry line. 
Therefore the polynomial $Q_{\n2+1}-Q_\n2=C_w$ must vanish on the legal 
branch of the curve. In general, there is no further factorisation,
as shown by the example $w=a^3$, for which the polynomial $C_w=a-1$ 
is irreducible.
If $n$ is even, then the points $z_{\n2+1}$ and $z_\n2$ are placed 
symmetrically with respect to the symmetry line, and hence 
$Q_{\n2+1}=Q_{\n2-1}$, that is, $C_w$ vanishes. 
For $w=a^4$ we have $C_w=a^2-2$, which is irreducible.
From theorem \ref{thm:Branches}, we see that 
$\n2=\lfloor n/2\rfloor$ branches of the curve $\cQ_n$ 
belong to the curve $\cC_w$, while the remaining $n-1-\n2$ branches
belong to ${\mathcal{\Cbar}}_w=0$.
The former comprises all parameters corresponding to paths of 
odd rank, while the latter those of even rank. 
Here the term rank refers to the \textit{orbital rank},
namely the number of half-turns around the origin, which,
as already noted, may not be related to the number of factors 
$ab$ and $ba$ in the word.

The existence of non-legal branches cannot be avoided by 
considering only irreducible curves.
For instance, the legal branch for the word $w=a^\kappa$
is the line $a=\zeta_\kappa$ [cf.~(\ref{eq:Segments})],
and $\zeta_\kappa$ is a root of the irreducible polynomial 
$\Psi_{2\kappa}(a)$ [cf.~(\ref{eq:Psi})]. 
This polynomial has degree $\phi(2\kappa)/2$, where $\phi$ is 
Euler's function \cite[p 37]{Niven}). For $\kappa\geqslant 4$ 
such a degree is greater than one, corresponding to as many 
branches; so there is a non-legal branch.

\subsection{Congruences}%
  \label{section:Functions}

In this section we consider functions defined on a critical curve 
$\cC_w$, with word $w=w_0\cdots w_{n-1}$. 
Unless indicated otherwise, the results of this section will apply 
to the more general setting of reduced palindromic words, as in 
proposition \ref{prop:Palindrome}.  
To lighten up the notation, we omit explicit reference to $w$ and write
$\cC$ for $\cC_w$ etc.

Given the polynomial $C(a,b)$ of a curve [see (\ref{eq:C})], 
we consider the polynomial ideal $\idl{C}=C(a,b)\mathbb{Z}[a,b]$ 
of all the multiples of $C$ in $\mathbb{Z}[a,b]$. 
(For background, see, e.g., \cite{CoxEtAl}.)
The quotient ring $\mathbb{Z}[a,b]/\idl{C}$ of residue classes 
modulo $\idl{C}$, namely the sets of the form $P+\idl{C}$ for 
$P\in\mathbb{Z}[a,b]$, represents the polynomial functions 
$P:\cC\to\mathbb{R}$.
We write $P\equiv Q\mod{C}$ to mean that $P-Q\in \idl{C}$, in which
case $P$ and $Q$ represent the same function on $\cC$.

Thus equation (\ref{eq:QSymmetry}) yields
\begin{equation}\label{eq:Symmetry}
Q_{n-t}(w_1,\ldots,w_{n-t-1})\equiv Q_t(w_1,\ldots,w_{t-1}) \mod{C},\qquad t=0,\ldots,n,
\end{equation} 
and considering that $Q_n \equiv Q_0=0$ and $Q_{n-1}\equiv Q_1=1$, on the 
curve $\cC$ the matrix (\ref{eq:Mws}) takes the form
\begin{equation}\label{eq:MonC}
\rM_{[w]} \equiv 
\begin{pmatrix}
  -Q_{n-1}(w_2,\ldots,w_{n-1}) & 0\\
  w_0-Q_{n-2}(w_2,\ldots,w_{n-2}) & -1
\end{pmatrix}
\equiv
\begin{pmatrix}
  -1 & 0\\
  w_0-Q_{n-2} & -1
\end{pmatrix}
\mod{C}
\end{equation}
where the second congruence follows from the fact that $\mathrm{det}(\rM)$ 
has unit determinant (cf.~\cite[theorem 3.4]{LagariasRainsI}). 

The variation of a function $f:\mathbb{R}^2\to \mathbb{R}$ 
along the curve $C$ is given by the Poisson brackets 
\begin{equation*}
\{f, C\}= \frac{\partial f}{\partial a}\frac{\partial C}{\partial b}
-\frac{\partial f}{\partial b}\frac{\partial C}{\partial a}.
\end{equation*}
We consider the observables
\begin{equation}\label{eq:phi_t}
\varphi_t=\arctan_2(Q_{t-1},Q_{t}),\qquad t=0,\ldots,n,
\end{equation}
which represent the angle of the
rays\footnote{$\varphi_t$ lies in 
the interval $(-\pi,\pi]$.} in the orbit of the map $\rF$, that is, 
the points of the orbit segment of the circle map.
This follows from the relation $x_t=\mathrm{sign}(w)Q_t$, where
$(x_0,\ldots,x_n)$ is the one-dimensional orbit segment associated
to the curve. 

We have
\begin{equation}\label{eq:Poisson}
\{\varphi_t, C\}=\frac{Q_t^2}{Q_t^2+Q_{t-1}^2}\{Q_{t-1}/Q_t,C\}
= \frac{Q_t\{Q_{t-1},C\}-Q_{t-1}\{Q_t,C\}}{Q_t^2+Q_{t-1}^2}.
\end{equation}
Thus $\{\varphi_t,C\}$ has the form $\Delta_t/\Vert z_t\Vert^2$, where
\begin{equation}\label{eq:Delta_t}
\Delta_t=Q_t\{Q_{t-1},C\}-Q_{t-1}\{Q_t,C\}, \qquad t=0,\ldots,n,
\end{equation}
and $\Vert z_t\Vert^2=Q_t^2+Q_{t-1}^2$. The latter has no real 
roots, because any common root of $Q_t$ and $Q_{t-1}$ would be
common to all $Q_i$s, which is impossible since $Q_1$ has no roots.

In preparation for the next statement, we consider the polynomials
\begin{equation}\label{eq:Xi}
\Xi_i
=
\begin{cases}
Q_i^2\frac{\partial C}{\partial a}& \mbox{\/if}\qad  w_i=b\\
-Q_i^2\frac{\partial C}{\partial b}& \mbox{\/if}\qad  w_i=a,
\end{cases}
\qquad i=0,\ldots,n-1
\end{equation}
If $\underline w$ is a palindrome, we have $\Xi_{n-i}=\Xi_i$, 
for $i=1,\ldots,n-1$.

We now establish formulae for $\Delta_t$.

\begin{theorem}\label{thm:Poisson}
Let $\cC$ be a boundary curve, and let $\Delta$ and $\Xi$ be as above. 
The following holds:
\vspace*{-20pt}
\begin{enumerate}
\item [i)] \qad $\displaystyle 
\Delta_0=0,\qquad
\Delta_t=\sum_{i=0}^{t-1}\Xi_i,
\quad t=1,\ldots,n,\qquad \Delta_n\equiv 0\mod{\cC}$
\end{enumerate}
\vspace*{-10pt}
If, in addition, $\cC$ is a critical curve, we have ($\n2=\lfloor n/2\rfloor$)
\vspace*{-10pt}
\begin{enumerate}
\item [ii)] \qad $\displaystyle
\Delta_{n-t}\equiv -\Delta_{t+1}\mod{C},\qquad
\Delta_t\equiv -\frac{1}{2}\sum_{i=t}^{n-t}\Xi_i\mod{C},
\quad t=1,\ldots,\n2.$
\newline\vskip 1pt\qad If $n=2\n2+1$, then
$\displaystyle \Delta_{\n2+1}=\{Q_\n2,Q_{\n2+1}\}C
   \equiv 0\mod{C}$.
\item [iii)] \qad $\displaystyle 
\sum_{t=1}^n \{\varphi_t,C\}\equiv 0\mod{C}$.
\end{enumerate}
\end{theorem}

\proof i) Using linearity and Liebnitz rule for Poisson brackets, 
(\ref{eq:Delta_t}) becomes
\begin{eqnarray*}
\Delta_t&=&(w_{t-1}Q_{t-1}-Q_{t-2})\{Q_{t-1},C\}
      -Q_{t-1}\{w_{t-1}Q_{t-1}-Q_{t-2},C\}\\
&=&Q_{t-1}\{Q_{t-2},C\}-Q_{t-2}\{Q_{t-1},C\}-Q_{t-1}^2\{w_{t-1},C\}.
\end{eqnarray*}
This gives
\begin{equation}\label{eq:Qrecursion}
\Delta_{t}=\Delta_{t-1}+Q_{t-1}^2\{C,w_{t-1}\},
\end{equation}
with 
$$
\{C,w_i\}=\begin{cases}
\displaystyle {\partial C}/{\partial a} & \mbox{if}\qad w_i=b\\
\noalign{\vskip 5pt}
\displaystyle -{\partial C}/{\partial b} & \mbox{if}\qad w_i=a.
\end{cases}
$$
From (\ref{eq:Delta_t}) we have $\Delta_0=\Delta_1=0$ and 
iterating the above recursion we find, for $t>0$
$$
\Delta_t=\sum_{i=0}^{t-1} Q_i^2\{C,w_i\}
    =\sum_{i=0}^{t-1}\Xi_i,
$$
as desired.
To establish the congruence, we let $Q_n=C\Cbar$ [see (\ref{eq:C})
and following remark], and compute
\begin{eqnarray*}
\Delta_n&=&C\Cbar\{Q_{n-1},C\}-Q_{n-1}\{C\Cbar,C\}\\
&=&C\Cbar\{Q_{n-1},C\}-Q_{n-1}C\{\Cbar,C\}\equiv 0\mod{C}.
\end{eqnarray*}
We have proved i).

We prove ii).
Using (\ref{eq:Symmetry}) and (\ref{eq:Delta_t}), we obtain 
\begin{equation}\label{eq:DeltaSymmetry}
\Delta_{t+1}\equiv -\Delta_{n-t}\mod{C},
\qquad t=0,\ldots,n-1,
\end{equation}
which, for $n=2\n2+1$, yields $\Delta_{\n2+1}\equiv 0\mod{C}$.
More precisely,
\begin{eqnarray*}
\Delta_{\n2+1}&=&Q_{\n2+1}\{Q_\n2,Q_{\n2+1}-Q_\n2\}
          -Q_\n2\{Q_{\n2+1},Q_{\n2+1}-Q_\n2\}\\
&=&\{Q_\n2,Q_{\n2+1}\}(Q_{\n2+1}-
          Q_\n2)=\{Q_\n2,Q_{\n2+1}\}C.
\end{eqnarray*}

Iterating (\ref{eq:Qrecursion}) backward and using i) we obtain 
the formula 
\begin{equation}\label{eq:Delta_tFormula2}
\Delta_t 
  \equiv-\sum_{i=t}^{n-1}\Xi_i\mod{C},\quad t=1,\ldots,n-1.
\end{equation}
Since $\Xi_{n-i}=\Xi_i$, we find, for $t=1,\ldots,\n2$
the above sum becomes 
\begin{eqnarray*}
\sum_{i=t}^{n-1}\Xi_i
&=& 
\sum_{i=t}^{n-t}\Xi_i+\sum_{i=n-t+1}^{n-1}\Xi_i\\
&\equiv&\sum_{i=t}^{n-t}\Xi_i+\sum_{i=1}^{t-1}\Xi_i
\mod{C}.
\end{eqnarray*}
Since $0=Q_0\equiv Q_n\mod{C}$, the range of the rightmost sum may 
be extended to include $i=0$.
Then, substituting the above expression in (\ref{eq:Delta_tFormula2}), 
and adding the latter to formula i), we obtain
$$
\Delta_t\equiv 
   -\frac{1}{2}\sum_{i=t}^{n-t}\Xi_i\mod{C}\qquad
t=1,\ldots,\n2,
$$
which completes the proof of ii).

We prove iii).
If $n$ is odd, then from ii) we obtain $\Delta_{\n2+1}\equiv0\mod{C}$.
Keeping this in mind, we find
\begin{eqnarray*}
\sum_{t=1}^n \{\varphi_t,C\}
&\equiv&
\sum_{t=1}^{\lfloor n/2\rfloor}\left(\frac{\Delta_{t}}{Q_{t}^2+Q_{t-1}^2}+
\frac{\Delta_{n-t+1}}{Q_{n-t+1}^2+Q_{n-t}^2}\right)\\
&\equiv&\sum_{t=1}^{\lfloor n/2\rfloor}
   \frac{\Delta_{t}+\Delta_{n-t+1}}{Q_{t}^2+Q_{t-1}^2}\equiv0\mod{C}.
\end{eqnarray*}
This establishes iii) and the proof of the theorem is complete.
\endproof

Some remarks on theorem \ref{thm:Poisson} are in place.
Part ii) says that $\Delta_{\n2+1}$ vanishes identically on ${\cC}$ for 
$n=2\n2+1$. This is due to time-reversal symmetry: the point
$z_{\n2+1}$ of the orbit segment never leaves the symmetry axis.

The statement iii) says that along an orbit segment of
a critical curve, the sum of the angles is constant.
To find the value of the constant, we represent the points of the 
orbit of $\rF$ as complex numbers $z_t=x_t+\mathrm{i}x_{t-1}$, 
with $z_0=-\mathrm{sign}(w)\mathrm{i}$.
We seek the value of
$$
S=\sum_{t=1}^n\varphi_t=
\sum_{t=1}^{n}\mathrm{arg}(z_t)\qquad z_t=\rF^t(z_0),
$$
for an arbitrary point on $\cC$.
We rewrite this sum as
\begin{equation}\label{eq:SumOfAngles}
S=\frac{1}{2}\sum_{t=0}^{n-1}
   [\mathrm{arg}(z_{t+1})+\mathrm{arg}(z_{n-t})].
\end{equation}
For any non-zero complex number $z$ 
we have $\arg(z)+\mathrm{arg}(\rR(z))\equiv\frac{\pi}{2}\mod{2\pi}$.
From reversibility, (\ref{eq:SumOfAngles}) becomes
$$
S=\frac{1}{2}\sum_{t=0}^{n-1}\frac{\pi}{2}\mod{2\pi}
 \equiv\frac{\pi n}{4}\mod{\pi}.
$$
In the above formula the modulus may be removed.
The value of the sum can be shown to depend on the numbers 
of rays lying in the third quadrant, which is constant along
$C$ as long as the code doesn't change.

In the next section we shall examine further geometrical consequences 
of theorem \ref{thm:Poisson}.

\section{Intersections of curves}\label{section:Intersections}

A \textbf{double point} is a point of intersection of two distinct
critical curves. 
Many geometrical properties of a critical curve $\cC$ are determined 
by its intersections with other critical curves.
For instance, $\cC$ has a single legal branch [theorem 
\ref{thm:Branches} i)], but in general not all points of 
that branch are legal. It turns out that the legal part of the branch,
which we call the \textbf{legal arc}, is delimited by certain double points.

\begin{theorem} \label{thm:LegalEagle}
For  $(\a,\b) \in (-2,2) \times (-2,2)$,
the legal points on the single legal branch of the critical curve 
$\cC_w$ of \eqref{eq:CriticalCurve}  are not isolated. 
They form legal arcs that are delimited by double points.
\end{theorem}
\proof
It suffices to consider the case of positive $w$.
From theorem \ref{thm:Branches} and its proof, we can assume the existence of a parameter pair $(\a_0,\b_0)$ that is a non-isolated point of  $\cC_w$, for which there is critical orbit, from $(0,-1)$ to $(1,0)$, with given positive word $w$ whose reduced word is the palindrome \eqref{eq:block}.  
We seek to investigate whether there persists a critical orbit with the same word for parameters in an open neighbourhood of $(\a_0,\b_0)$ on $\cC_w$. 

From \eqref{eq:MwReduced}, we have for $z_1=(1,0)$ and $t>1$, 
$z_{t}=(Q_t(w_1,\ldots,w_{t-1}), Q_{t-1}(w_1,\ldots,w_{t-2}))$, 
depending on the first $t-1$ symbols of $\underline w$. 
The coordinates of $z(t)$ are polynomial, hence continuous, 
functions of $\a$ and $\b$. 
By assumption, when $\a=\a_0 \in (\zeta_{\kappa-1}, \zeta_{\kappa})$, 
$\kappa \geqslant 2$ and $\b=\b_0 \in (\zeta_{\ell-1}, \zeta_{\ell})$, 
$\ell \geqslant 2$, we have $z_{n}=(0,1)$. 
From proposition \ref{prop:legalA}, the domain $\bd(\rM_a)$ is divided into 
$\kappa$ sectors by $L^{\pm}$ and the $\kappa-1$ rays ${\bf r}_j(\a)$ 
of \eqref{eq:SecRay} that are functions of $\a$ and rotate anti-clockwise 
as $\a$ increases. Likewise the domain $\bd(\rM_b)$ is divided into $\ell$ 
sectors by $L^{\pm}$ and the analogous $\ell-1$ rays that depend on $\b$ 
only and  rotate anti-clockwise as $\b$ increases.  
The existence of a critical orbit with positive word can be viewed as 
the occurrence of a non-empty intersection of the forward semi-infinite 
$L^{-}$-orbit with the backward semi-infinite $L^{+}$ orbit.
Under the assumption that $t=n$ is the first visit to $(0,1)$, we have 
that the forward orbit of $(0,-1)$ avoids all sector boundaries in $\bd(\rM_a)$ and $\bd(\rM_b)$ until it coincides with the first ray $(U_{\kappa-1}(\a), U_{\kappa}(\a))$  in the final visit to $\bd(\rM_a)$ before completing the critical segment (otherwise it would arrive at $L^{-}$ or $L^{+}$ earlier than claimed).  
Because the sector boundaries in $\bd(\rM_a)$ or $\bd(\rM_b)$ are all 
iterates of their first ray, a critical (boundary) orbit for a positive word 
exists if and only if the forward orbit of $(0,-1)$ coincides at some point 
with the first ray of $\bd(\rM_a)$.
As soon as this happens, the forward orbit of $(0,-1)$ and the backward 
orbit of $(0,1)$ coincide along a finite critical orbit segment. 
Necessarily, this condition on parameters $(\a, \b)$ of coinciding 
with the first ray of $\bd(\rM_a)$ is equivalent to belonging to $\cC_w$. 
 
As the critical orbit segment is finite there are a finite number of possible first ray collisions that can occur as we vary $\a$ and $\b$. 
If at $(\a_0,\b_0)$, the only collision is with the first ray 
of $\bd(\rM_a)$ in the final block, we can maintain this collision 
by staying on $\cC_w$ and, by continuity of the orbit rays in $\a$ 
and $\b$, continue to avoid first ray collisions in other blocks of 
the word on some open arc containing $(\a_0, \b_0)$.  
This arc will be delimited by parameter values corresponding to other 
first ray collisions before the final block, i.e., double points. 
\endproof

We defer the question of whether there is a unique legal arc on the legal branch to a later paper \cite{RobertsEtAl}.  This has to do with whether the earlier first ray collisions are transverse. 
We claim each of the rays in the orbit $\rF^j \,(0,-1), j \geqslant 1$, move clockwise as $a$ or $b$ increases, i.e., the opposite direction to the first rays which move anti-clockwise.  To see this, realise that the analysis in theorem \ref{thm:Poisson} applies for arbitrary $\cC$, not just the critical curve.  Taking $\cC$ to be a horizontal or vertical line shows that $\{\varphi_t, b\}=\frac{\partial \varphi_t}{\partial a}$ and $\{\varphi_t, a\}=-\frac{\partial \varphi_t}{\partial b}$ are both negative.  So moving to the right and down on the critical curve is necessary to maintain the final ray at $(0,1)$.

Let $c=(\a,\b)\in\cC$ be such that for some $t\not=0,n$ we have 
$Q_t(c)=0$. 
We collect all values of $t$ for which 
$Q_t(c)=0$ to form the finite sequence 
\begin{equation}\label{eq:T}
T=T_w(c)=(t_1,t_2,\ldots)\qquad t_j<t_{j+1},
\end{equation}
called the \textbf{intersection sequence} of the curve at $c$.
In preparation for the next statement, we 
denote by $Q_{t,u}$ the $t$-th $Q$-polynomials (\ref{eq:PQ}) 
for the word $u$, and by $\cQ_{t,u}$ the corresponding curve.
As before, we write $Q_t$ for $Q_{t,w}$.

\begin{theorem}\label{thm:DoublePoints}
Let a critical curve $\cC_w$ have non-empty intersection 
sequence $T(c)=(t_1,t_2,\ldots)$ at a point $c$.
Then $c$ is a double point, $|T|$ is even, and
the orbit at $c$ is periodic with minimal period $t_2$.
If we let, for $j=1,\ldots,|T|/2$,
\begin{equation}\label{eq:uv}
w=uvu',\quad\mbox{where}\quad u=w_0\cdots w_{t_j-1}
  \quad\mbox{and}\quad|u'|=|u|,
\end{equation}
then $c$ lies at the intersection of $\cC_w$ and three curves, 
namely 
\begin{equation}\label{eq:Cuv}
\cC_{u,j}=\cQ_{t_j,u},
\qquad
\cC_{v,j}=\cQ_{n-2t_j,v},
\qquad
\cC_{uv,j}=\cQ_{n-t_j,{uv}}.
\end{equation}
The orbital rank of $\cC_{v,j}$ is odd, while those of $\cC_{u,j}$ 
and $\cC_{uv,j}$ have, respectively, the same and the opposite 
parity as $j$.
\end{theorem}

\proof Let $t\in T$. Then $n-t\in T$, from (\ref{eq:Symmetry}).
However, we cannot have $t=n-t=n/2$, because then from $Q_{n/2}=0$ 
and $C=Q_{n/2+1}-Q_{n/2-1}=0$ we would have that all $Q_i$s vanish 
at $c$, but $Q_1=1$ does not. So $|T|$ is even.
By symmetry ($\underline w$ is a palindrome), the $t$-th 
and the $(n-t)$-th rays are distinct boundary rays, so at $c$ 
the boundary segment visits both rays twice, and is therefore 
periodic. Moreover one of $\cQ_t$ and $\cQ_{n-t}$ is a critical
curve, and hence $c$ is a double point.
During one period the orbit must visit both boundary rays, 
so the period cannot be $t_1$, which is minimal. 
For the same reason, the rays visited at consecutive 
$t_j$s must be different. 
It follows that the minimal period of the orbit is $t_2$,
and that the orbital rank of $\cQ_{t_j}$ has the same
parity as $j$, while that of $\cQ_{n-t_j}$ has opposite 
parity. 

Consider now the decomposition (\ref{eq:uv}) for fixed $j$. 
Since both $u$ and $uv$ are prefixes of $w$, with $|u|=t_j$
and $|uv|=n-t_j$, we find that $\cC_{u,j}=\cQ_{t_j}$ and 
$\cC_{uv,j}=\cQ_{n-t_j}$. This establishes the parity of
the ranks of $\cC_{u,j}$ and $\cC_{uv,j}$.
To show that at $c$ the orbit on the curve $\cC_{v,j}$
is the same as the middle segment of the orbit on the 
boundary curve $\cC$, we must verify that the
polynomial $Q_{n-2t_j,v}$ has the correct initial 
conditions prescribed by (\ref{eq:PQ}).

There are two cases. 
If $j$ is odd, then $\cC_{u,j}$ is a critical curve.
From (\ref{eq:Symmetry}) we then have 
$Q_{1,u}(c)=Q_{t_j-1,u}(c)=1$, and hence
$(Q_{t_j}(c),Q_{{t_j}-1}(c))=(0,1)=(Q_{0,v}(c),-Q_{-1,v}(c))$.
Therefore, for $t=t_j,\ldots,n-t_j$ we have $Q_{t}(c)=-Q_{t-t_j,v}(c)$, 
and in particular $Q_{n-2t_j,v}(c)=Q_{n-t_j}(c)=0$.
Thus $\cC_{v,j}$ is a boundary curve of odd orbital rank.

If $j$ is even, then $\cC_{u,j}$ is not a critical curve,
and to compute $Q_{t_j-1}(c)=Q_{-1,v}(c)$ we cannot use (\ref{eq:Symmetry}).
However, since $t_j$ is a (not necessarily minimal) period, 
for some $0<t<t_j$ the $t$-th point of the orbit must visit 
the end ray of $\cC$. Then, by concatenating
two critical curves [equivalently, by composing two matrices
of type (\ref{eq:MonC})], we conclude that $Q_{-1,v}=\pm 1$.
The analysis proceeds as before, and we conclude again that 
$\cC_{v,j}$ is a boundary curve with odd orbital rank.
\endproof

Theorem \ref{thm:DoublePoints} identifies a prominent
set of double points associated with a critical curve $\cC_w$.
They occur at intersections with curves of lower rank,
corresponding to factors of $w$.
Several points are worth considering.

\noindent
i) Since $Q_t(c)=Q_{n-t}(c)=0$, the symbols $w_t$ and $w_{n-t}$ 
may be changed independently without affecting the dynamics. 
As a result, at a legal double point the code $w$ may be,
and typically is, \textit{improper} 
(see the beginning of section \ref{section:CriticalCurves}).
\newline
ii) The $Q$-polynomials of odd rank may be replaced by 
the corresponding $C$-divisor, according to proposition 
\ref{prop:Palindrome} i), eliminating the branches with even
orbital rank.
\newline
iii) From the palindrome property of $w$ and proposition 
\ref{prop:ContinuedFractions} ii), we find that the words
$u$ and $u'$ generate the same curve, as so do
$uv$ and $vu'$. So the sub-words $u,v$, and $uv$ 
describe completely the decomposition (\ref{eq:uv}).
\newline
iv) At a double point $c$ of a critical curve $\cC_w$
there are $|T_w(c)|/2$ distinct decompositions, each
involving intersections of boundary curves of lower rank.
Thus the minimum number of decomposition is 1, while the maximum
is $(\iota-1)/2$, where $\iota$ is the rank of $w$.

For illustration, consider the rank 7 palindrome let $w=(a^3b^4)^3a^2$
[cf.~theorem \ref{thm:FirstGenerationCurves}, ii), section 
\ref{section:RotationalDomains}]. 
At the double point $c=(0,2\cos(\pi/5))$, we have $T=(2,7,9,14,16,21)$,
so the period is equal to 7. 
The symbols $w_i$ are improper for $i=2,9,16$, so 
the proper code is $(a^2b^5)^3a^2$
[see theorem \ref{thm:FirstGenerationCurves} ii) 1 below].
Theorem \ref{thm:DoublePoints}, applied to the proper code,
gives $|T|/2=3$ decompositions:
\begin{equation}\label{eq:Decompositions}
\begin{array}{ccccccc}
t&&u&v&u'&\quad&\mbox{ranks}\\
\noalign{\vskip 1pt\hrule\vskip 3pt}
2&&a^2&(b^5a^2)^2b^5&a^2&&\mbox{1,5,1}\\
7&&a^2b^5&a^2b^5a^2&b^5a^2&&\mbox{2,3,2}\\
9&&a^2b^5a^2&b^5&a^2b^5a^2&&\mbox{3,1,3}
\end{array}
\end{equation}
The proper code is a prefix of the periodic 
word $(a^2b^5)^\infty$ of period 7.

At $c=(1,2\cos(3\pi/11))$ we have $T(c)=(3,20)$, so the period is 20.
The symbol $w_{20}=w_{n-3}$ is improper. The proper code is
$(a^3b^4)^2a^3b^3a^3$, which is a prefix of the periodic word
$((a^3b^4)^2a^3b^3)^\infty$ of period $20$.
We have $|T|/2=1$ decomposition:
$$
t=3,\quad u=a^3,\quad v=(b^4a^3)^2b^3,\quad u'=a^3 \qquad\mbox{ranks:}\, 1,5,1.
$$

Next we provide a formula for the rotation number at a double point,
and a partial converse of theorem \ref{thm:DoublePoints}.

\begin{lemma} \label{lma:Intersection}
Let $\cC_w$ and $\cC_{w'}$ be distinct critical curves 
of ranks $\iota_w$ and $\iota_{w'}$ (with $\iota_{w'}\leqslant \iota_w$, say), 
which intersect at a legal point c. 
Then
\begin{equation}\label{eq:OppoSign}
\theta(c)=
\begin{cases}
\displaystyle \frac{\iota_w+\iota_{w'}}{2(|w|+|w'|)}& \mbox{if \/} \mathrm{sign}(w)\not=\mathrm{sign}(w')\\
\displaystyle \frac{\iota_w-\iota_{w''}}{2(|w|-|w''|)}& \mbox{if \/} \mathrm{sign}(w)=\mathrm{sign}(w'),
\end{cases}
\end{equation}
where $w''=w'$ if $\iota_{w'}\not=\iota_w$, and otherwise $w''$ is 
any prefix of $w'$ whose length is an odd-order element of the
intersection sequence of $w'$ at $c$. 
\end{lemma}

\proof
Let $w$ and $w'$ have opposite sign, with $w$ positive (say). 
Then the orbit segment of $w$ maps $L^-$ to $L^+$ in 
$|w|$ iterates, and that of $w'$ maps $L^+$ to $L^-$ in $|w'|$ 
iterates. Thus at the double point $c$ the (non necessarily
minimal) period is $|w|+|w'|$. 
Let $\iota_w$ and $\iota_{w'}$ be the rank of the words $w$ 
and $w'$ at $c$. 
Since $c$ is legal, the rotation number at $c$ is equal to half the combined 
rank divided by the period, regardless of whether the codes are 
proper or improper, the ranks being unaffected by this property.

Let now $w$ and $w'$ have the same sign (positive, say), hence the same initial
ray at $c$.
We have two cases. If $\iota_{w'}<\iota_w$, then the word $w'$ is equivalent 
to a prefix of $w$, so $w\sim w'u'$ for some non-empty word $u'$, and $L^+$ 
is periodic under $u'$, with (not necessarily minimal) period $|u'|=|w|-|w'|$, 
and even rank $\iota_{u'}=\iota_w-\iota_{w'}$.
Letting $w''=w'$ and computing the rotation number from period and rank, 
as above, gives the desired formulae.

If $\iota_{w'}=\iota_{w}$, then since the two orbit segments 
have the same initial condition, we have $|w|=|w'|$. Because
$w\not=w'$, the two words will differ at some boundary ray, 
that is, their common intersection sequence $T$ at $c$ 
is non-empty. In particular, $T$ has some odd-order element $t_j$.
It now suffices to let $w''$ be the prefix of $w'$ of length $t_j$ 
and proceed as above. 
\endproof

The intersection of critical curves described in 
lemma \ref{lma:Intersection} always leads to
a decomposition $w=uvu'$ of type (\ref{eq:uv}). 
Indeed, in the equal sign case, letting $n=|w|$ and $n-t=|w''|$
we have $Q_{n-t,w}(c)=0$, and, by symmetry, $Q_{t,w}(c)=0$. 
Then both $t$ and $n-t$ belong to the intersection sequence $T_w(c)$. 
The required decomposition is obtain by letting
$u=w_0\cdots w_{t'-1}$ where $t'=\min(t,n-t)$.
In the unequal sign case we also obtain a decomposition of type 
(\ref{eq:uv}), by considering the word $ww'w$.

To study intersections of curves, we introduce a sequence of 
vectors in $\Z[a,b]^2$:
\begin{equation}\label{eq:Gamma}
\cG=(\Gamma_0,\ldots,\Gamma_{n-1})
\end{equation}
where
\begin{equation}\label{eq:Gamma_t}
\Gamma_t=\sum_{i=0}^t\gamma_i,\qquad
\gamma_t=\begin{cases} (Q_t^2,0) & \mbox{if}\qad w_t=a\\
                       (0,Q_t^2) & \mbox{if}\qad w_t=b.
         \end{cases}
\end{equation}
A key property of this sequence is derived from theorem \ref{thm:Poisson}.
From i) we have $\Delta_n\equiv 0\mod{C}$, and 
using again i) and (\ref{eq:Xi}) we find
$$
\sum_{i=1}^{n-1} \Xi_i=\frac{\partial C}{\partial a}
\sum_{\stackrel{0\leqslant i<n}{w_i=b}} Q_i^2
-\frac{\partial C}{\partial b}
\sum_{\stackrel{0\leqslant i<n}{w_i=a}} Q_i^2\equiv 0\mod{C}.
$$
The rightmost congruence expresses the vanishing of a determinant 
on $\cC$, which establishes the following geometrical fact.

\begin{corollary}\label{cor:Parallel}
If $\cC$ is a boundary curve, then the vectors
\begin{equation}\label{eq:Parallel}
\nabla C
 =\Bigl(\frac{\partial C}{\partial a},\frac{\partial C}{\partial b}\Bigr)
\qquad\mbox{and}\qquad
\Gamma_{n-1}
\end{equation}
are parallel at every point of the curve.
\end{corollary}
Thus the normal to a boundary curve may be determined without 
computing derivatives. More precisely, there is a rational 
function $\lambda\in\Q(a,b)$, such that $\nabla C=\lambda\Gamma_{n-1}$ 
on $\cC$. Since both vectors in (\ref{eq:Parallel}) are non-zero, 
the function $\lambda$ is regular and non-zero on $\cC$.
From the corollary we also deduce at once that if the rank of $\cC$ 
is greater than one, then the partial derivatives of $C$ have the same sign 
[cf.~theorem \ref{thm:Branches} ii)].

\begin{figure}[tb]
\centering
\begin{tikzpicture}[scale=1]

\begin{axis}[
   xmin=0, xmax=9,
   ymin=0, ymax=11,
   xtick={0,2,4,6,8},
   ytick={0,2,4,6,8,10},
   xticklabels={0,2,4,6,8},
   yticklabels={0,2,4,6,8,10},
   axis lines=middle,
   x axis line style=->,
   y axis line style=->,
]
\addplot[very thick,color=black,mark=*] coordinates{
(0., 0.)
(1., 0.)
(2.008, 0.)
(2.008, 0.6724e-4)
(2.008, .9869)
(2.008, 2.692)
(2.008, 3.201)
(2.149, 3.201)
(3.335, 3.201)
(3.851, 3.201)
(3.851, 3.336)
(3.851, 4.772)
(3.851, 6.209)
(3.851, 6.344)
(4.368, 6.344)
(5.554, 6.344)
(5.695, 6.344)
(5.695, 6.852)
(5.695, 8.558)
(5.695, 9.545)
(5.695, 9.545)
(6.703, 9.545)
(7.703, 9.545)
};
\addplot[thick,color=red] coordinates{
(0., 0.)
(7.703, 9.545)
};
\addplot[thick,color=blue] coordinates{
(2.008, 0.)
(5.695, 9.545)
};
\end{axis}
\end{tikzpicture}
\caption{\label{fig:Polygonal}\rm\small
The odd-rank legal polygonal $\cG_w(c)$ for the word $w=(a^3b^4)^3a^2$ at
the intersection $c=(1,2\cos(3\pi/11))$ of the curve $\cC_w$ 
with $\cC_{a^3}$. The added line segments are the medians of $\cG_w$ and 
of $\cG_v$ with $v=b^4(a^3b^3)^2$, which represent the normal to the 
curves $\cC_w$ and $\cC_v$ at $c$. 
These words and curves are described in theorems \ref{thm:DoublePoints}
and \ref{thm:FirstGenerationCurves} ii).
}
\end{figure}
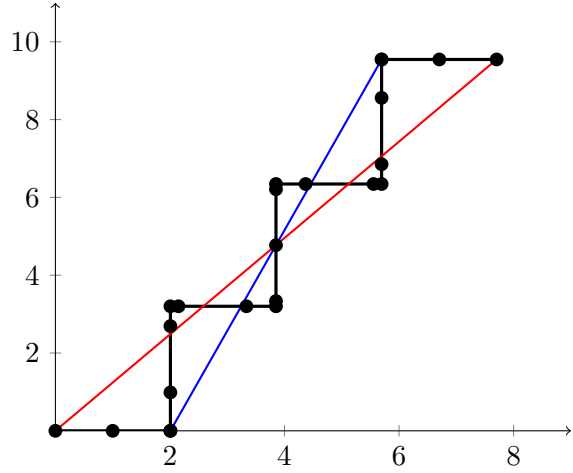

For every parameter pair $c=(\a,\b)\in\R^2$, the sequence (\ref{eq:Gamma})
defines a polygonal on the plane (still denoted by $\cG$),
obtained by connecting the elements of the sequence by line 
segments (figure \ref{fig:Polygonal}).
If $c\in\cC_w$, then $\cG(c)$ is called the \textbf{polygonal of $\cC_w$}
(or of $w$). The terms `legal' for points on curves, and `rank' for 
words or curves, will also be used for polygonals.

By construction [see (\ref{eq:Gamma})] the vertices of $\cG$ 
correspond to code changes, that is, to the occurrence of 
the factor $ab$ or $ba$ in $w$. Therefore $\cG$ has $\iota$ 
edges (line segments) and $\iota-1$ vertices, where $\iota$ is the rank of $w$.
If $\cC$ is a critical curve, then from (\ref{eq:Symmetry}) and 
the fact that the reduced word is a palindrome, we have 
$\gamma_i\equiv \gamma_{n-i}\mod{C}$.
As a consequence, the polygonal of a critical curve is 
symmetrical with respect to its barycentre.

The elements of $\cG$ are not necessarily distinct, and if
$\Gamma_{t}=\Gamma_{t-1}$ for some $t$, then we say that $\Gamma_t$
is an \textbf{intersection point} of $\cG$.
If $\Gamma_t$ is an intersection point of $\cG$, then 
$\gamma_t=(0,0)$, whence $Q_t(c)=0$. 
Since $1=Q_1\equiv Q_{n-1}\mod{C}$, the end-points of $\cG$ are 
not intersection points, so $Q_t$ lies in the interior of the 
boundary segment, that is, $c$ is a double point of $\cC$, and
$t$ belongs to the intersection sequence $T(c)$, see (\ref{eq:T}). 
This argument may be reversed, to show that at a double point, 
the intersection points of the polygonal and the elements of 
the intersection sequence are in bi-unique correspondence.

Intersection points may occur at an arbitrary position on the polygonal. 
However, a legal intersection point necessarily corresponds to a 
code change, even if the code is improper. 
Thus the intersection points of a legal polygonal $\cG$ must
occur at the vertices, and they occur in pairs, symmetrically 
placed with respect to the barycentre of $\cG$.
Thus all vertices of a legal polygonal are intersection points 
if and only if the intersection sequence is maximal: 
$|T|=\iota-1$. 

Let $\cG$ be a polygonal of odd rank. The segment joining its
end-points will be called the \textbf{median} of $\cG$, which bisects
$\cG$'s middle segment. $\cG$ is said to be \textbf{regular} if 
no intersection point of $\cG$ lies on the median. 
To formulate a sufficient condition for regularity, we consider the
intersection sequence $T_w(\c)=(t_1,t_2,\ldots)$ [cf.~(\ref{eq:T})]
of a curve $\cC_w$ at a point $\c$.
We say that $T_w(\c)$ is \textbf{simple} if the
the word $w_0\cdots w_{t_1-1}$ is equivalent to a word of rank 1
at $\c$.
An empty intersection sequence will also be considered simple.

\begin{theorem}\label{thm:Regular}
An odd-rank legal polygonal with simple intersection sequence is regular. 
\end{theorem}
\proof 
First we show that the intersection points of an odd-rank polygonal 
either all lie on the median, or none of them does.
Let $T=(t_1,t_2,\ldots)$ be the intersection sequence 
of $\cG=(\Gamma_0,\ldots,\Gamma_{n-1})$.
The statement is trivially true if $T$ is empty (e.g., for 
rank-one polygonals), since there are no intersection points.
If $T$ is non-empty, then the orbit at $c$ is periodic with 
minimal period $t_2$, by theorem \ref{thm:DoublePoints}.
By periodicity, we have
$\Gamma_{t_{2k}}=k\Gamma_{t_2}$, $k=1,\ldots,|T|/2$, that is,
all even-rank intersection points lie on the segment joining
the origin $\Gamma_0$ to the last intersection point $\Gamma_{t_{|T|}}$.
By symmetry, all odd-rank intersection points lie on the segment
joining $\Gamma_{n-1}$ to the first odd-rank intersection point
$\Gamma_{t_1}$. Therefore $\Gamma_{t_1}$ lies on the median,
if and only if all intersection points of $\cG$ lie on the median, 
as desired. (Note that this statement holds even if the polygonal 
is not legal.)

It now suffices to show that that if $\cG$ is legal and $T$ is 
simple, then one vertex of $\cG$ does not belong to the median.
For $\iota>1$ (the rank-1 case holds by definition)
$\cG$ has at least two vertices, and all the intersection 
points (if any) are at the vertices.
Since $T$ is simple, the first vertex is an intersection point,
which not on the median by construction.
The proof of the theorem is complete.
\endproof

Next we show that regular polygonals bring about
transversal intersections of curves.

\begin{lemma}\label{lma:Transversality}
With the notation of theorem \ref{thm:DoublePoints},
let the polygonal of a critical curve $\cC_w$ be 
regular at a legal double point $c$.
Then the following holds:
\vspace*{-10pt}
\begin{enumerate}
\item [i)] $\cC_w$ intersect transversally the curves
$\cC_{u,j},\cC_{v,j},\cC_{uv,j}$, for $j=1,\ldots,|T|/2$.
\item [ii)] All pairwise intersections of the families of curves 
\begin{equation*}
ii.1)\qad \cC_{u,j},\cC_{v,j},\cC_{uv,j},\qad\mbox{for fixed  $j$};
\hskip 40pt
ii.2)\qad \cC_{v,j},\cC_{v,j'},\qad\mbox{for any $j\not=j'$}.
\end{equation*}
are transversal.
\item [iii)] If $|T|>2$, then the curves $C_{uv,j},C_{u,k}$ 
with $j$ odd and $k$ even are tangent.
\end{enumerate}
\end{lemma}
\proof
We fix $j$ in (\ref{eq:uv}), and let $k(j)=|T|-j+1$. 
Then $j$ and $k$ have opposite parity, while $\Gamma_{t_j}$ and 
$\Gamma_{t_k}$ are symmetric with respect to the centre of $\cG$. 
We define three subsequences of $\cG$:
\begin{equation}\label{eq:GuGv}
\cG_{u,j}=(\Gamma_0,\ldots,\Gamma_{t_j-1}),\quad
\cG_{v,j}=(\Gamma_{t_j},\ldots,\Gamma_{t_k-1}),\quad
\cG_{uv,j}=(\Gamma_0,\ldots,\Gamma_{t_k-1}).
\end{equation}
We want to show that $\cG_{u,j},\cG_{v,j}-\Gamma_{t_j}$ and $\cG_{uv,j}$
are, respectively, the polygonals of the curves $\cC_{u,j},
\cC_{v,j}$ and $\cC_{uv,j}$; by construction these polygonals 
are then embedded in $\cG$.

For $u$ and $uv$ the result is immediate, since these words 
are prefixes of $w$.
For the word $v$, we must verify that
for $i=-1,\ldots,n-2t_j-1$ the sequences $Q_{t_j+i}(c)$ and
$Q_{i,v}(c)$ agree in absolute value. Indeed the 
corresponding codes are the same, and the verification
that the initial conditions ($i=-1,0$) are the same
has already been given in the proof of theorem \ref{thm:DoublePoints}.
It follows that $\cG_{v,j}$ is an embedding of the polygonal
of $\cC_{v,j}$ into $\cG$ by translation (see figure \ref{fig:Polygonal}).
 
Next we deal with transversality.
From corollary \ref{cor:Parallel}, the intersections of two curves 
is transversal if and only if the medians of the corresponding 
polygonals are not parallel; so we shift our attention to medians.

The medians of the polygonals $\cG_{u,j},\cG_{v,j}$
and $\cG_{uv,j}$ are given by
$$
A_{u,j}=[\Gamma_0,\Gamma_{t_j-1}]
\qquad
A_{v,j}=[\Gamma_{t_j},\Gamma_{t_k-1}]
\qquad
A_{uv,j}=[\Gamma_0,\Gamma_{t_k-1}].
$$
By symmetry, the $A_{v,j}$s have a common mid-point.
Since $\cG$ is regular, from the proof of theorem 
\ref{thm:Regular} we have that one end-point of 
$A_{v,j}$ belongs to the segment joining $\Gamma_0$ to 
$\Gamma_{t_{|T|}}$, while the other belongs 
to the segment joining $\Gamma_{n-1}$ to $\Gamma_{t_1}$.
For this reason, no two $A_{v,j}$ can be parallel, which is ii.2).

Let $A=[\Gamma_0,\Gamma_{n-1}]$ be the median of $\cG$.
Since $\cG$ is regular, $A=[\Gamma_0,\Gamma_{n-1}]$ and 
$A_{v,j}$ intersect transversally at their common mid-point. 
Thus $A$ and $A_{v,j}$ are the diagonals of a parallelogram 
which has one vertex at the origin $\Gamma_0$. 
The medians $A_{u,j}$ and $A_{uv,j}$ of the other polygonals
connect the origin to the vertices of $A_{v,j}$, so no pair of
medians can be parallel. This is ii.1).

In total, we obtain $|T|/2$ distinct non-degenerate parallelograms, 
sharing one diagonal $A$, with their second diagonals $A_{v,j}$ 
forming a pencil through the centre of the parallelogram.
This suffices to establish that $A$ is transversal to 
all medians $A_{u,j}, A_{v,j}$ and $A_{uv,j}$, which implies
the transversal intersections of the corresponding curves. This is i).

Suppose that $|T|>2$. Then the medians
$$
A_{uv,1},A_{u,2},A_{uv,3},A_{u,4},\ldots
$$
have one end-point at $\Gamma_0$, while (as pointed out earlier) 
the other end-points are collinear.
This implies the tangential intersection of the corresponding 
curves, which is iii).

The proof of the lemma is complete.
\endproof

\begin{lemma}\label{lma:Equivalence}
Let $\cC$ be a critical curve, and let $c\in\cC$ be a legal point. 
The following statements are equivalent:
\vspace*{-10pt}
\begin{enumerate}
\item [i)] For some $t$ and $s\not=t$, we have  $\varphi_t(c)-\varphi_s(c)=0$ 
[cf.~(\ref{eq:phi_t})].
\item [ii)] The intersection sequence $T(c)$ is non-empty.
\end{enumerate}
\end{lemma}
\proof If $\varphi_t(c)-\varphi_s(c)=0$, with $t<s$, say, then, 
from the continuity and invertibility of the circle map, we have 
that $\varphi_{t+j}(c)-\varphi_{s+j}(c)=0$, $j=-t,\ldots,n-s$.
Thus the orbit segment is part of a periodic orbit of period 
$k:=s-t$, which contains both $\varphi_0$ and $\varphi_n$, that is, $c$ is 
a double point.
The value $j=-t$ corresponds to a collision of the $k$th ray with the initial
boundary ray ($z_k=z_0$), while $j=n-s$ yields the same phenomenon for the end 
boundary ray ($z_{n-k}=z_n$). Thus $Q_k(c)=0$, which means that $k\in T(c)$,
that is, $T(c)$ is non-empty.
Conversely, if $T(c)$ is non-empty, then by theorem \ref{thm:DoublePoints}
we have $\varphi_{t_2}(c)-\varphi_0(c)=0$.
\endproof

Condition i) expresses the collision of points of the
orbit of the circle map at $c$. The transversality ---or
lack of it--- of such a collision is expressed by the validity
---or lack of it--- of the inequality $\varphi_t^\prime(c)\not=\varphi_s^\prime(c)$, 
where $\varphi^\prime=\{\varphi,C\}$.
This condition is independent from the choice of $t$ and $s$ along the orbit.
Indeed, suppose that $\varphi_t^\prime(c)-\varphi_s^\prime(c)=0$.
Then (\ref{eq:Poisson}) and (\ref{eq:Delta_t}) give 
$(Q_t^2+Q_{t-1}^2)\Delta_s=(Q_s^2+Q_{s-1}^2)\Delta_t$. 
Using the recursion (\ref{eq:Qrecursion}), and keeping in mind 
that $w_s=w_t$ and that $Q_sQ_{t-1}=Q_tQ_{s-1}$, we obtain 
$\Delta_{s\pm 1}(Q_t^2+Q_{t-1}^2)=\Delta_{t\pm 1}(Q_s^2+Q_{s-1}^2)$. 
Since $z_t=\lambda z_s$, for some $\lambda$, then also 
$z_{t\pm 1}=\lambda z_{s\pm 1}$, from the local linearity of $\rF^{\pm 1}$. 
Thus $\Delta_{s\pm 1} \Vert z_{t\pm 1}\Vert^2 =\Delta_{t \pm 1} \Vert z_{s\pm 1}\Vert^2$, that is,
$\{\varphi_{s\pm 1},C\}=\{\varphi_{t\pm 1},C\}$.
Repeating this argument an appropriate number of times,
we find that all intersections are non-transversal. 

An \textbf{end-point} of a curve is a point on the boundary of 
the curve's legal arc.

\begin{theorem} \label{thm:Endpoints}
At an end-point $c$ of a critical curve $\cC$ the intersection
sequence $T(c)$ is non-empty. Conversely, if $T(c)$ is non-empty and
$\cG(c)$ is regular, then $c$ is an end-point. \end{theorem}
\proof
The initial and final rays of a critical curve are boundary rays, which,
by definition, remain fixed along the curve. 
If $c\in\cC$ is an end-point of $C_w$, then at $c$ some 
intermediate ray must become a boundary ray. 
Equivalently, the intersection sequence $T_w(c)$ is non-empty.

Suppose that $T(c)$ is non-empty. Then, according to theorem
\ref{thm:DoublePoints}, the curve $\cC$ intersects the curves of
lower rank given in (\ref{eq:Cuv}).
If $\cG(c)$ is regular, then these curves intersect $\cC$ 
transversally, from lemma \ref{lma:Transversality} i).
This in turn means that the corresponding rays intersect boundary rays
transversally [that is, $\{\varphi,C\}\not=0$ at $c$, see (\ref{eq:Poisson})], 
so that the code becomes illegal at $c$.
\endproof

The regularity of $\cG$ is not necessary for $c$ to be an
end-point. Indeed using lemma \ref{lma:Equivalence}, the end-points 
may be characterised in terms of collisions of points of the 
orbit of the circle map, and these collisions need not be transversal
to render a code illegal.

\section{First-generation critical curves}\label{section:RotationalDomains}

Proposition \ref{prop:Segments} established all critical curves of
rank 1. In this section we construct two classes of boundary curves 
of higher rank. 
First, all rank-2 curves, obtained by concatenating two rank-1 
words of opposite sign. 
Second, all the critical curves of the \textbf{first generation}.
They are constructed by concatenating an arbitrary number of copies 
the same rank-2 word, and then extending the resulting word 
---to the right or to the left--- with a suitable rank-1 word.
First-generation curves are regrouped to form infinite pencils,
incident to the same point, the \textbf{basis} of the pencil.
For reason of brevity, we shall only deal with curves which lie in the
first quadrant of the $(\a,\b)$-parameter space, leaving the general
case to the sequel of this paper \cite{RobertsEtAl}.

\begin{figure}[tb]
\centering
\begin{tikzpicture}[scale=5.0]
\tikzmath{\c2=0; \c3=1; \c4=1.414213562; \c5=1.618033989; \c6=1.732050808;}
\tikzmath{\mydot=0.015; \d=0.04; \dd=0.08;}
\fill [gray!60] (\c3,\c4) rectangle (\c4,\c5);

\draw[black, thin] (\c2,\c3) -- (\c5,\c3) -- (\c5,\c6) -- (\c2,\c6) -- (\c2,\c3);
\draw[black, thin] (\c2,\c4) -- (\c5,\c4);
\draw[black, thin] (\c2,\c5) -- (\c5,\c5);
\draw[black, thin] (\c2,\c6) -- (\c5,\c6);
\draw[black, thin] (\c3,\c3) -- (\c3,\c6);
\draw[black, thin] (\c4,\c3) -- (\c4,\c6);
\draw[black, thin] (\c5,\c3) -- (\c5,\c6);

\draw [black, very thick] (\c3,\c5) to [bend left=15] (\c4,\c4);
\draw [red, very thick] (\c2,\c5) .. controls (1.09,1.55)..  (\c4,\c4);
\draw [red, very thick] (\c3,\c5) .. controls (1.35,1.56) and (1.5,1.55) .. (\c5,\c4);
\draw [blue, very thick] (\c3,\c6) to [bend left=20] (\c4,\c4);
\draw [blue, very thick] (\c3,\c5) to [bend left=20] (\c4,\c3);

\filldraw[black] (\c4,\c4) circle (\mydot);
\filldraw[black] (\c3,\c5) circle (\mydot);
\filldraw[black] (\c2,\c5) circle (\mydot);
\filldraw[black] (\c4,\c3) circle (\mydot);
\filldraw[black] (\c5,\c4) circle (\mydot);
\filldraw[black] (\c3,\c6) circle (\mydot);
\filldraw[black] (\c4,\c4+0.125) circle (\mydot*0.6);
\filldraw[black] (\c3+0.233,\c5) circle (\mydot*0.6);
\filldraw[black] (\c3+0.291,\c5-0.048) circle (\mydot*0.5);
\filldraw[black] (\c3,\c4+0.125) circle (\mydot*0.6);
\filldraw[black] (\c3+0.120,\c4+0.105) circle (\mydot*0.5);
\filldraw[black] (\c3+0.220,\c4) circle (\mydot*0.6);

\node[black] at (\c4,\c3-\d) {\tiny $(4,3)$};
\node[black] at (\c3,\c6+\d) {\tiny $(3,6)$};
\node[black] at (\c2-\dd,\c5) {\tiny $(2,5)$};
\node[black] at (\c3-\dd,\c5+\d) {\tiny $(3,5)$};
\node[black] at (\c4+\dd,\c4-\d) {\tiny $(4,4)$};
\node[black] at (\c5+\dd,\c4) {\tiny $(5,4)$};

\end{tikzpicture}
\caption{\label{fig:CriticalCurves}\rm\small
Illustration of theorem \ref{thm:FirstGenerationCurves} for the domain 
$\cD_{3,4}$ (the grey rectangle).
The black curve, of type i) in the theorem, is the axis of the domain.
The blue and red curves are critical curves, of type ii) and iii), 
respectively and $m=1$; they intersect the boundary of the domain 
at the mediant rotation number of the adjacent vertices.
The integer pairs represent the subscripts of the vertices $\bv_{\kappa,\ell}$,
which are joined by the curves. 
}
\end{figure}
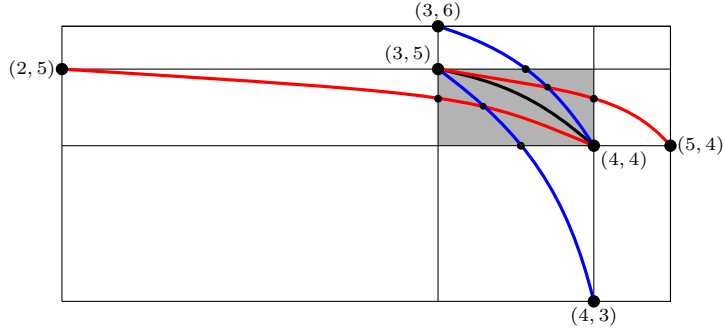

We begin by partitioning the first quadrant into
rectangular \textbf{rotational domains} $\cD_{\kappa,\ell}$,
given by $\zeta_\kappa \leqslant \a\leqslant \zeta_{\kappa+1}$
and $\zeta_\ell \leqslant \b\leqslant \zeta_{\ell+1}$
(figure \ref{fig:RotDomains}).
We will show that the legal arc of a
rank-2 word lies within a single rotational domain, while 
that of a first-generation critical curve occupies two neighbouring 
domains (see figure \ref{fig:CriticalCurves}).
 
We establish some notation.
The boundary words $a^\kappa$ and $b^\ell$\/ have rank 1 
and opposite sign. Thus, from equation (\ref{eq:Segments}) and
lemma \ref{lma:Intersection}, at the intersection of the corresponding
curves we have the double point
\begin{equation}\label{eq:ckl}
\bv_{\kappa,\ell}:=\bC_{a^\kappa}\cap \bC_{b^\ell}=(\zeta_\kappa,\zeta_\ell)
\hskip 40pt
\theta(\bv_{\kappa,\ell})=\frac{1}{\kappa+\ell}.
\end{equation}
We shall make repeated use of the functions
\begin{equation}\label{eq:ThetaZeta}
\vartheta^\pm_{j,m}=\frac{m}{mj\pm1},
\hskip 30pt
\zeta^\pm_{j,m}=2\cos \pi\vartheta^\pm_{j,m}.
\end{equation}

In the rest of this paper we write $\cC_w$ to mean the legal branch 
of $\cC_w$. 

We now state and prove the main result of this section.

\goodbreak
\begin{theorem}\label{thm:FirstGenerationCurves}
~~

\vspace*{-5pt}
\begin{enumerate}
\item [i)] For all $\kappa,\ell\geqslant 2$ the words 
$a^{\kappa+1}b^{\ell}$ and $b^{\ell+1}a^{\kappa}$ 
give the same rank-2 curve with 
end-points $\bv_{\kappa,\ell+1}$ and $\bv_{\kappa+1,\ell}$.
\vspace*{5pt}
\item [ii)] The following sequences of words $w^{(m)}, m=1,2,\ldots$
define pencils of critical curves, with end-points $\be$ 
(the basis) and $\be_m$, given by
\begin{equation*}
\begin{array}{lcccccc}
&&w^{(m)}&\be&\be_m&\theta(\be_m)&\\
\noalign{\vskip 6pt}
1)&&
(a^{\kappa+1}b^{\ell})^ma^\kappa
&
\bv_{\kappa,\ell+1}
&
(\zeta_{\kappa+1},\zeta^-_{\ell,m})
&
\vartheta^-_{\kappa+\ell+1,m}
&
\raisebox{-12pt}{\quad$\kappa\geqslant 2, \ell>2$}\\
\noalign{\vskip -6pt}
2)&&
a^{\kappa+1}(b^{\ell}a^{\kappa})^m
&
\bv_{\kappa+1,\ell-1} 
&
(\zeta_{\kappa},\zeta^+_{\ell,m})
&
\vartheta^+_{\kappa+\ell,m}
&
\\
\noalign{\vskip 7.5pt}
3)&&
(b^{\ell+1}a^{\kappa})^mb^\ell
&
\bv_{\kappa+1,\ell}
&
(\zeta^-_{\kappa,m},\zeta_{\ell+1})
&
\vartheta^-_{\kappa+\ell+1,m}
&
\raisebox{-12pt}{\quad$\kappa> 2, \ell\geqslant 2$,}\\
\noalign{\vskip -6pt}
4)&&
b^{\ell+1}(a^{\kappa}b^{\ell})^m
&
\bv_{\kappa-1,\ell+1}
&
(\zeta^+_{\kappa,m},\zeta_\ell)
&
\vartheta^+_{\kappa+\ell,m}
&
\end{array}
\end{equation*}
while $\theta(\be)$ is given in (\ref{eq:ckl}).
\end{enumerate}
\end{theorem} 

\proof We denote by $\Sigma_i$, $i=1,\ldots,4$ the $i$th quadrant 
in phase space. 

We prove i).
The reduced words of $w=a^{\kappa+1}b^{\ell}$ and 
$w'=b^{\ell+1}a^{\kappa}$ are mapped into one another by a reflection
symmetry. From proposition \ref{prop:ContinuedFractions} ii), it then
follows that the corresponding curves are the same.

We consider the double points $\bv_{\kappa,\ell+1}$ and 
$\bv_{\kappa+1,\ell}$, with proper positive codes 
$(a^{\kappa} b^{\ell+1})^\infty$ and 
$(a^{\kappa+1} b^{\ell})^\infty$, respectively. 
At the former point, $w\sim a^{\kappa} b^{\ell+1}$;
at the latter, $w$ is proper.
Likewise, $w'$ is proper at $\bc_{\kappa,\ell+1}$ 
and $w'\sim b^\ell a^{\kappa+1}$.
Thus both points are legal for both codes and from theorem 
\ref{thm:Branches} i), there is precisely one  arc of $\cC_w$
connecting them.
As we proceed along this arc from $\bv_{\kappa,\ell+1}$ to 
$\bv_{\kappa+1,\ell}$, the $\kappa$-th ray 
$(Q_\kappa,Q_{\kappa-1})$ of the positive orbit
rotates clockwise into $\Sigma_1$, since $\kappa>1$.
The $(\kappa+1)$-st ray remains in $\Sigma_2$, because, by 
construction, there are 
no active branches of positive rank-1 curves in the interior of 
$\cD_{\kappa,\ell}$. 

Thus the legal arc of $\cC_w$ contains $\cC_w\cap \cD_{\kappa,\ell}$.
The monotonicity of ray rotations with $\a$ prevents the prolongation of the 
legal arc outside $\cD_{\kappa,\ell}$.
An analogous argument shows that the legal arc of $\cC_{w'}$ is the 
same as that of $\cC_{w}$. Statement i) is proved.

Next we turn to critical curves, making a preliminary remark. 
If a positive curve $\cC_w$ of rank greater than 1 intersects 
$\cC_{a^\kappa}$ at a legal point $\c$, 
then $\c$ is necessarily an end-point of $\cC_w$.
Indeed, since the code at $\c$ is $a^\kappa b\cdots$, 
the word $a^\kappa$ is equivalent to a prefix of $w$, and hence 
at $\c$ the intersection sequence of $\cC_w$ is simple.
Then, by theorem \ref{thm:Regular} the polygonal 
of the curve is regular at $\c$, and hence $\c$ is an 
end-point for $\cC_w$, from theorem \ref{thm:Endpoints}.
The same holds for a negative curve intersecting $\cC_{b^\ell}$.

We prove ii) part 1. Let $w^{(m)}=(a^{\kappa+1}b^{\ell})^m a^\kappa$.
At the double point $\bc_{\kappa,\ell+1}=(\zeta_\kappa,\zeta_{\ell+1})$,
the proper code is $(a^\kappa b^{\ell+1})^\infty$,
and $w^{(m)}\sim (a^\kappa b^{\ell+1})^m a^\kappa$. 
Thus, for all $m$, the point $\bc_{\kappa,\ell+1}$ is a legal point of the 
curve $\cC_{w^{(m)}}$, hence an end-point, from the above remark.
We have shown that $\be=\bc_{\kappa,\ell+1}$.

We now prolong $\cC_{w^{(m)}}$ inside the domain 
$\cD_{\kappa,\ell-1}\cup\cD_{\kappa,\ell}$. We write
$\bC_{\kappa,\ell}$ for the axis $\cC_{a^{\kappa+1}b^\ell}$ and we let
\begin{equation}\label{eq:L}
L:=\rM_{[w^{(m)}]} L^-=\rM_{[a^\kappa]}M_{[(a^{\kappa+1}b^{\ell})^m]}L^-.
\end{equation}
We fix $\a$ in the range $\zeta_\kappa<\a<\zeta_{\kappa+1}$,
and we proceed by induction on $m$.

Let $m=1$. We choose $\b$ so that $(\a,\b)$ lies in the region bounded 
by $\bC_{\kappa,\ell}$ and $\bC_{\kappa,\ell-1}$.  
If $(\a,\b)\in\bC_{\kappa,\ell}$, then i) and (\ref{eq:L}) give
$L=\rM_a^\kappa L^-\in\Sigma_1$, with the code $w^{(1)}$.
If $(\a,\b)\in\bC_{\kappa,\ell-1}$, then
\begin{equation}\label{eq:onLowerCurve}
\rM_{[a^{\kappa+1}b^{\ell-1}]}L^-=L^-
\qquad\mbox{and}\qquad
\rM_{[a^{\kappa+1}]}L^-\in \Sigma_2.
\end{equation}
An infinitesimal increase in $\b$ changes the code 
$a^{\kappa+1}b^{\ell-1}a^{\kappa+1}$ to $w^{(1)}$,
while leaving the orbit segment unchanged.
From (\ref{eq:onLowerCurve}) we conclude that just 
above $\bC_{\kappa,\ell-1}$ we have $L\in\Sigma_2$.
Thus, from continuity and monotonicity, there must 
be a unique value of $\b$ for which $L=L^+$, 
which shows that $w^{(1)}$ is a boundary word of rank $3$,
and that the curve $C_{w^{(1)}}$ has a legal branch between
$\bC_{\kappa,\ell-1}$ and $\bC_{\kappa,\ell}$.

Suppose now that for some $m\geqslant 1$ the curve $\cC_{w^{(m)}}$ 
has a legal branch between $\bC_{\kappa,\ell-1}$ and $\bC_{\kappa,\ell}$.
We fix $\a$ and choose $\b$ so that $(\a,\b)$ lies
between $\cC_{w^{(m)}}$ and $\bC_{\kappa,\ell}$.
If $(\a,\b)\in\bC_{\kappa,\ell}$, then i) and (\ref{eq:L}) give
$\rM_{[w^{(m+1)}]}L^-=\rM_{[a^\kappa]}L^-\in\Sigma_1$.
If $(\a,\b)\in \cC_{w^{(m)}}$ then, by the inductive hypothesis 
$\rM_{[w^{(m)}]} L^-=L^+$, while 
$\rM_{[b^{\ell+1}a^\kappa]}L^+\in\Sigma_2$. 
An infinitesimal increase in $\b$ causes 
the code $w^{(m)}b^{\ell+1}a^\kappa$ to become $w^{(m+1)}$,
without changing the orbit segment; thus on $\cC_{w^{(m)}}$
we have $\rM_{[w^{(m+1)}]}L^-\in \Sigma_2$.
It then follows that there is a unique value of $\b$
in the specified range for which $\rM_{[w^{(m+1)}]}L^-=L^+$,
which completes the induction.

The point $\c_m=\cC_{w^{(m)}}\cap\cC_{a^{\kappa+1}}$
has the form $(\zeta_{\kappa+1},\b)$, for some 
$\b=2\cos 2\pi\vartheta_b$ to be determined.
Such a point is legal by continuity, hence is an
end-point of $\cC_{w^{(m)}}$ (by the remark above),
so that $\be_m=\c_m$.
To determine $\b$ we first find $\theta(\c_m)$.
The words $w^{(m)}$ and $a^{\kappa+1}$ have the same sign, and 
$w^{(m)}-a^{\kappa+1}=b^\ell(a^{\kappa+1}b^\ell)^{m-1}a^\kappa$;
lemma \ref{lma:Intersection} then gives
\begin{equation}\label{eq:Thetan}
\theta(\c_m)=\frac{m}{m(\kappa+\ell+1)-1}=\vartheta^-_{\kappa+\ell+1,m}.
\end{equation}
The parameter $\vartheta_\b$ is related to the rotation number 
$\theta(\zeta_{\kappa+1},\b)$ by equation (\ref{eq:ThetaRhoSegments}). 
Solving the latter for $\vartheta_\b$ and 
using (\ref{eq:Thetan}), we obtain $\b=\zeta^-_{\ell,m}$, 
that is $\be_m=(\zeta_{\kappa+1},\zeta^-_{\ell,m})$, as desired.

We have shown that for any $m\geqslant1$ the curve $\cC_{w^{(m)}}$ has a
legal arc in $\cD_{\kappa,\ell-1}\cup\cD_{\kappa,\ell}$ 
with the stated end-points. The proof of ii) 1 is complete.

We prove ii) part 2. Let $w^{(m)}=a^{\kappa+1}(b^\ell a^\kappa)^m$. 
An argument analogous to that used in ii) 1 shows that
$\bv_{\kappa+1,\ell-1}$ is a legal point of the curve 
$\cC_{w^{(m)}}$ for all $m$, so that $\be=\bv_{\kappa+1,\ell-1}$.

To prolong $\cC_{w^{(m)}}$ inside the domain 
$\cD_{\kappa,\ell}\cup\cD_{\kappa,\ell-1}$ we proceed by induction
on $m$, keeping in mind that the base case is already established,
since for $m=1$ the words for the two statements are the same. 
Thus assume that for some $m\geqslant 1$ the curve $\cC_{w^{(m)}}$ 
has a legal branch between $\bC_{\kappa,\ell-1}$ and $\bC_{\kappa,\ell}$.
We fix $\a$ in the range $\zeta_{\kappa}<\a<\zeta_{\kappa+1}$ and 
$\b$ so that $(\a,\b)$ lies between $\bC_{\kappa,\ell-1}$ and $\cC_{w^{(m)}}$.

If $(\a,\b)\in\bC_{\kappa,\ell-1}$, then from part i)
$\rM_{[(b^\ell a^{\kappa})^{m+1}]}$ maps $L^+$ to itself,
and hence $\Sigma_2$ to itself. 
The second equation in (\ref{eq:onLowerCurve}) then 
implies that $\rM_{[w^{(m+1)}]}L^-\in\Sigma_2$.
If $(\a,\b)\in\cC_{w^{(m)}}$, then, by the inductive hypothesis 
$\rM_{[w^{(m)}]} L^-=L^+$, while $\rM_{[b^{\ell}a^\kappa]}L^+\in\Sigma_1$. 
Since $w^{(m)}b^\ell a^\kappa=w^{(m+1)}$, on $\cC_{w^{(m)}}$ we have 
$\rM_{[w^{(m+1)}]}L^-\in\Sigma_1$.
It then follows that there is a unique value of $\b$
in the specified range for which $\rM_{[w^{(m+1)}]}L^-=L^+$,
which completes the induction.

As above, $\be_m$ is an end-point of the curve, and the formulae for $\be_m$ 
and $\theta(\be_m)$ are computed as the corresponding formulae in ii) 1.
The proof of ii) 2 is complete.

The proof of ii) 3,4 is obtained from that of ii) 1,2 by merely 
exchanging parameters.
\endproof

Let us examine theorem \ref{thm:FirstGenerationCurves} at the light of 
the material of section \ref{section:Intersections}.
We consider the curves of case ii) 1, of rank $2m+1$ 
the other cases being analogous.
At the left end-point 
$\be=(\zeta_\kappa,\zeta_{\ell+1})$ we have
$w=(a^{\kappa+1}b^\ell)^ma^\kappa\sim(a^\kappa b^{\ell+1})^ma^\kappa$.
With reference to the decomposition (\ref{eq:uv}), the $Q$-polynomials
of the $2m$ prefixes 
$u=a^\kappa,a^\kappa b^{\ell+1},a^\kappa b^{\ell+1} a^\kappa,\ldots,
(a^\kappa b^{\ell+1})^m$ vanish, leading to the following 
maximal intersection sequence
$$
T(\be)=(\kappa,\kappa+\ell+1,2\kappa+\ell+1,2\kappa+2(\ell+1),
    \ldots,m\kappa+m(\ell+1)).
$$
From theorem \ref{thm:Regular} and \ref{thm:Endpoints} we conclude
that the polygonal $\cG(\be)$ is regular and that all these 
intersections are transversal.

At the right-end point we find, from (\ref{eq:ThetaZeta}) and 
ii) 1:
$$
\be_m=(\zeta_{\kappa+1},\zeta^-_{\ell,m}),
\qquad \zeta^-_{\ell,m}=2\cos(\pi m/(m\ell-1)).
$$
We now show that
$$
T(\be_m)=(t_1,t_2)=(\kappa+1,m-\kappa-1).
$$
The factor $a^{\kappa+1}$ produces a rotation by $\pi$,
so $t_1=\kappa+1$. The factor $b^\ell$ produces a
rotation by an angle greater than $\pi$. For $m>1$, the quantity
$mt/(m\ell-1)$ is an integer for $t=m\ell-1=(m-1)\ell+\ell-1$, 
and no smaller positive $t$, as easily verified. It follows that
$T(\be_m)$ has no term $t$ between $t_1$ and $t_2=m-t_1$, and that
at $\be_m$ we have 
$w\sim(a^{\kappa+1}b^\ell)^{m-1}a^{\kappa+1}b^{\ell-1}a^{\kappa+1}$.
Thus the intersection sequence is maximal only for $m=1$.
However, the polygonal $\cG(\be_m)$ is still regular for every $m$,
because the first intersection point lies on the first vertex.

When combined with lemma \ref{lma:Intersection}, theorem 
\ref{thm:FirstGenerationCurves}
has the following immediate corollary, which gives the rotation number 
of the points of intersection of the odd-rank curves in parts ii) 
within the domain $\cD_{\kappa,\ell}$. 
We extend the parameters to include the case $m=0$,
which corresponds to the words $a^\kappa$ and $b^\ell$, respectively
(see proposition \ref{prop:Segments}).
This is legitimate, since we only make use of the fact that
all these curves are legal in $\cD_{\kappa,\ell}$.

\begin{corollary}\label{cor:Grid} 
For any $\kappa,\ell\geqslant 2$ and $n,m\geqslant 0$, 
the following holds
\begin{equation}
\begin{array}{ccc}
w&w'&\theta(\cC_w\cap \cC_{w'})\\
\noalign{\vskip 5pt\hrule\vskip 5pt}
(a^{\kappa+1} b^\ell)^na^\kappa & (b^{\ell+1} a^\kappa)^m b^\ell
     &\theta^-_{\kappa+\ell+1,n+m+1}\hfill\\
\noalign{\vskip 3pt}
a^{\kappa+1} (b^\ell a^\kappa)^n & b^{\ell+1}(a^{\kappa+1} b^\ell)^m
    &\theta^+_{\kappa+\ell+1,n+m+1}\hfill\\
\end{array}
\end{equation}
\end{corollary}
Note that the rotation number depends on the parameters through
the sums $\kappa+\ell$ and $n+m$.

\section*{Acknowledgements}
JAGR and AS kindly thank Queen Mary, University of London, for their hospitality.  JAGR and FV thank Valerie Berth\'e for stimulating discussions on rotational words. 
JAGR thanks Tim Siu for his help in creating the parameter space plot in Figure 1.
This research was supported by the Australian Research Council and 
by JSPS KAKENHI Grant No.~JP16KK0005.

\section*{Appendix: polynomial identities}

For the purpose of factoring various polynomials appearing in
our analysis, we introduce the sequence of polynomials
\begin{equation}\label{eq:Psi}
\Psi_1(X)=X-2,\quad \Psi_2(X)=X+2,\hskip 30pt
\Psi_n(x+x^{-1})=C_n(x)x^{-\phi(n)/2}\qquad n=3,4,\ldots,
\end{equation}
where $C_n(X)$ is the $n$-th cyclotomic polynomial
(that is, the roots of $C_n$ are the primitive $n$-th roots of unity 
\cite[section 2.4]{McCarthy}),
and $\phi$ is Euler's function \cite[p 37]{Niven}.
For $n>2$, $\Psi_n$ is a monic polynomial in $X=x+x^{-1}$, of degree $\phi(n)/2$. 
Moreover, $\Psi_n$ is irreducible for all $n$, and its roots are the 
distinct numbers $2\cos(2\pi k/n)$, with $k$ coprime to $n$.
These properties of $\Psi_n$ are established from the fact that 
the polynomial $C_n$ has degree $\phi(n)$, is irreducible and
reflexive\footnote{meaning that $X^{\phi(n)}\rF_n(X^{-1})=\rF_n(X)$},
together with the repeated use of the identity
\begin{equation}\label{eq:xkIdentity}
x^k+x^{-k}=(x+x^{-1})(x^{k-1}+x^{1-k})-(x^{k-2}+x^{2-k})\qquad k\in\Z.
\end{equation}

\subsection{Iterates}\label{section:Iterates}
Let now $\rM_X$ be as in (\ref{eq:M}).
A straightforward induction shows that
the matrix $\rM_X^n$ can be written as
\begin{equation}\label{eq:M^n}
\rM_X^n =
\begin{pmatrix}
  \cp_{n+1} & -\cp_{n}\\
  \cp_{n} & -\cp_{n-1}
\end{pmatrix},\qquad n\geqslant 0,
\end{equation}
where $\cp_n$ satisfies the recursion relation
\begin{equation}\label{eq:Chebyshev}
\cp_{-1}(X)=-1, \quad \cp_{0}(X)=0, \hskip 40pt 
\cp_{n+1}(X)=X \cp_{n}(X) -\cp_{n-1}(X), \quad n\geqslant 0.
\end{equation}
We see that for $n\geqslant 1$, $\cp_n$ is a polynomial in $X$ with 
integer coefficients and degree $n-1$; the term of degree 
$k$ is nonzero if and only if $k$ has the same parity as $n-1$.
The recursion (\ref{eq:Chebyshev}) is a special case of the more general relation
$$
\cp_{n}=\cp_{k}\cp_{n-k+1}-\cp_{k-1}\cp_{n-k}\qquad n,k\in\Z
$$
which is obtained from (\ref{eq:M^n}) and the identity $\rM^n=\rM^{n-k}\rM^k$.
Using (\ref{eq:Chebyshev}), one sees that for $n>0$, $\cp_{n}(X)={\bar U}_{n-1}(X/2)$, 
where $\bar U_n$ is the $n$th Chebyshev polynomial of the second kind, 
whence $\cp_n(2\cos \theta)=\frac{\sin n\theta}{\sin \theta}$.


\end{document}